\documentclass[journal]{IEEEtran}
\usepackage{amssymb}
\usepackage{amsmath}
\usepackage{amscd}
\usepackage{graphicx}
\usepackage{graphics}
\usepackage{xcolor,graphicx}
\usepackage{csquotes}
\usepackage{mathtools}
\usepackage{caption}
\usepackage{float}
\usepackage{subcaption}
\ifCLASSINFOpdf
\else
\fi
\begin{document}
%
\title{Compound Regularization of Full-waveform Inversion for Imaging Piecewise Media}
%
%
%

\author{Hossein S. Aghamiry,~\IEEEmembership{}
        Ali Gholami,~\IEEEmembership{}
        St\'ephane Operto,~\IEEEmembership{}

\thanks{Hossein S. Aghamiry is with the Institute of Geophysics, University of Tehran, Tehran, Iran, and University Cote d'Azur - CNRS - IRD - OCA, Geoazur, Valbonne, France,
  e-mail: h.aghamiry@ut.ac.ir, aghamiry@geoazur.unice.fr. }
  \thanks{Ali Gholami is with the Institute of Geophysics, University of Tehran, Tehran, Iran,
  e-mail: agholami@ut.ac.ir.}
\thanks{St\'ephane Operto is with University Cote d'Azur - CNRS - IRD - OCA, Geoazur, Valbonne, France,
  e-mail: operto@geoazur.unice.fr.}
         }

%
%


\markboth{A preprint, March~2019}%
{Shell \MakeLowercase{\textit{et al.}}: Bare Demo of IEEEtran.cls for Journals}

%



\maketitle

\begin{abstract}
The nonlinear and ill-posed nature of full waveform inversion (FWI) requires us to use sophisticated regularization techniques to solve it. In most applications, the model parameters may be described by physical properties (e.g., wave speeds, density, attenuation, anisotropic parameters) which are piecewise functions of space. Compound regularizations are thus necessary to reconstruct properly such parameters by FWI. We consider different implementations of compound regularizations in the wavefield reconstruction inversion (WRI) method, a formulation of FWI that extends its search space and prevent the so-called cycle skipping pathology. Our hybrid regularizations rely on Tikhonov and total variation (TV) functionals, from which we build two classes of hybrid regularizers: the first class is simply obtained by a convex combination (CC) of the two functionals, while the second relies on their  infimal convolution (IC).  In the former class, the model of parameters is required to simultaneously satisfy different priors, while in the latter the model is broken into  its basic components, each satisfying a distinct prior (e.g. smooth, piecewise constant, piecewise linear). We implement these types of compound regularizations in the WRI optimization problem using the alternating direction method of multipliers (ADMM). Then, we assess our regularized WRI in the framework of seismic imaging applications. Using a wide range of subsurface models, we conclude that compound regularizer based on IC leads to the lowest error in the parameter reconstruction compared to that obtained with the CC counterpart and the Tikhonov and TV regularizers when used independently.
\end{abstract}



\section{Introduction}

Full waveform inversion (FWI) seeks to estimate constitutive parameters of a medium with a wavelength-scale resolution by minimization of a distance between the recorded and simulated data subject to the wave equation is satisfied. This technology was originally developed in geophysical imaging \cite{Tarantola_1984_ISR}, and has spread more recently in other fields of imaging sciences such as medical imaging \cite{Sandhu_2015_FDU} or oceanography \cite{Wood_2008_FWI}.
FWI is a partial-differential equation (PDE) - constrained nonlinear optimization problem, which is classically solved with local reduced-space optimization methods \cite{Haber_2000_OTS}. In this linearized framework, a challenging source of non linearity is the so-called cycle skipping pathology. When FWI is solved with a classical least-squares difference-based distance \cite{Virieux_2009_OFW,Virieux_2017_FWI}, cycle skipping drives FWI to spurious local minima as soon as the initial parameter model does not allow to match the data with a kinematic error smaller than half a period. Other sources of error are noise, approximate wave physics and ill-posedness resulting from parameter cross-talk, coarse acquisition sampling and uneven illumination of the targeted structure. Designing regularization techniques that mitigate these source of errors is therefore a key challenge for the success of FWI applications.

A proper regularization should be driven by the shape and statistical characteristics of the medium to be imaged.
For example, in geophysical imaging, the subsurface can be represented by a piece-wise smooth medium, that is a model which contains smoothly varying and blocky components. The widespread Tikhonov regularizations \cite{Tikhonov_2013_NMF} rely on the smoothness assumption and hence fail to recover sharp interfaces of such media. Conversely, total variation (TV) regularizations \cite{Rudin_1992_NTV} are based on blockiness assumption and hence are more suitable to image large contrasts. However, they generate undesirable staircase imprint in smooth regions.
Regions characterized by smoothly-varying properties and those containing sharp contrasts have indeed different statistical properties. The former are characterized by the normal prior, while the latter by a heavy tailed prior. Consequently, simultaneous  recovery of both properties is difficult when one type of regularization is used (Tikhonov, TV, etc). 
To overcome this issue, a combination of different regularizations can be used \cite{Gholami_2010_ROL,Benning_2018_MRM}. 
A naive approach consists of the simple additive coupling or convex combinations (CC) of regularizations.  Alternatively, \cite{Gholami_2013_BCT} proposed to decompose explicitly the model into several components of different statistical properties and use an appropriate regularization to reconstruct each component. Using this strategy, they combined Tikhonov and TV regularizations (referred to as TT regularization) to reconstruct piece-wise smooth media. The smooth components are captured by the Tikhonov regularization, while the blocky ones are determined by the TV counterpart. In many applications, it has been shown that such a compound regularizations based upon infimal convolution (IC) outperforms the one based upon additive coupling  \cite{Bergmann_2018_PCF}. \\
TT regularization based upon IC has been successfully applied to FWI for seismic subsurface imaging in the framework of iteratively-refined wavefield reconstruction inversion (IR-WRI) \cite{Aghamiry_2018_HTT}. IR-WRI extends the search space of FWI and decreases cycle skipping accordingly by relaxation of the wave-equation constraint \cite{VanLeeuwen_2013_MLM,Aghamiry_2019_IWR,Aghamiry_2018_IFW}. Taking advantage of the bilinearity of the wave equation, IR-WRI breaks down FWI into two linear subproblems which are solved in an alternating mode: wavefield reconstruction driven by the observables and model-parameter estimation by minimization of the source residuals the relaxation generated. 
The linearity of the parameter-estimation subproblem provides a suitable framework to implement sophisticated nonsmooth regularizations. \\
In this study, in the following of \cite{Gholami_2013_BCT} and \cite{Aghamiry_2018_HTT}, we develop a general framework to combine a couple of regularization terms in IR-WRI through IC. 
Then, we specifically develop this framework for Tikhonov and TV regularizations, which are suitable for seismic subsurface imaging applications. Compared to \cite{Aghamiry_2018_HTT}, we jointly update the blocky and the smooth components through a variable projection process rather than in alternating mode. We first show that our new IC-based TT regularization outperforms the CC-based counterpart with several well-documented numerical benchmarks in the field of seismic imaging. We also compare the results obtained with these two TT regularizations with those obtained with total generalized variation (TGV) regularization, a combination of first and second order TV, and those obtained with Tikhonov and TV regularizations when used independently.
\section{Notation}
The mathematical symbols adopted in this paper are as follows.
We use italics for scalar quantities, boldface lowercase letters for vectors, and boldface capital letters for matrices and tensors.
We use the superscript $T$ to denote the adjoint of an operator. 
The $i$th component of the column vector $\bold{x}$ is shown by $\bold{x}_i$ and its absolute value is returned by $|\bold{x}_i|$.
For the real-valued $n$-length column vectors $\bold{x}$ and $\bold{y}$ 
the dot product is defined by $\langle \bold{x},\bold{y}\rangle=\bold{x}^T\bold{y}=\sum_{i=1}^n\bold{x}_i\bold{y}_i$ and
their Hadamard product, denoted by $\bold{x}\circ \bold{y}$, is another vector made up of their component-wise products, i.e. $(\bold{x}\circ \bold{y})_i=\bold{x}_i\bold{y}_i$.
The $\ell_1$- and $\ell_2$-norms of $\bold{x}$ are, respectively, defined by $\|\bold{x}\|_2=\sqrt{\langle \bold{x},\bold{x}\rangle}=\sqrt{\sum_{i=1}^n\bold{x}_i^2}$ and
$\|\bold{x}\|_1=\sum_{i=1}^n|\bold{x}_i|$.
For a real-valued $n\times m$ matrix $\bold{A}$, we use the (entrywise) $\ell_1$- and $\ell_2$-norms similar to vector norms as 
$\|\bold{A}\|_2=\sqrt{\sum_{i=1}^n\sum_{j=1}^m\bold{A}_{ij}^2}$ and $\|\bold{A}\|_1={\sum_{i=1}^n\sum_{j=1}^m|\bold{A}_{ij}|}$.

\section{Method}
In this section, we briefly review the frequency-domain FWI as a bi-convex feasibility problem and describe the reduced and  extended forms of FWI. We show how the latter can be solved with the alternating direction method of multipliers (ADMM) \cite{Boyd_2011_DOS} for a general regularization functional. 
\subsection{Full-waveform inversion}
Frequency-domain FWI with a general regularization term and bounding constraints can be formulated as \cite{Aghamiry_2019_IWR,Aghamiry_2019_IBCb}
\begin{equation} \label{main}
\begin{split}
&\underset{\bold{u},\bold{m}\in \mathcal{C}}{\min}  ~~~~~ \Phi(\bold{m})  \\
&\text{subject to} ~~ \bold{A(m)u}=\bold{b},~~ \bold{Pu}= \bold{d},
\end{split}
\end{equation} 
where $\bold{m} \in \mathbb{R}^{n\times 1}$ gathers unknown squared slowness, $\Phi(\bold{m})$ is an appropriate regularization term which we assume to be convex, $\mathcal{C} = \{\bold{x} \in \mathbb{R}^{n\times 1}~\vert~ \bold{m}_{{l}} \leq \bold{x} \leq \bold{m}_{{u}}\}$ is the set of all feasible models bounded by the lower bound $\bold{m}_{{l}}$ and the upper bound $\bold{m}_{{u}}$.\\ The first constraint in \eqref{main},
$\bold{A(m)u}=\bold{b}$, is a partial-differential equation (PDE) wherein $\bold{u} \in \mathbb{C}^{n\times 1}$ is the  wavefield and $\bold{b} \in \mathbb{C}^{n\times 1}$ is the source term. In this study, $\bold{A(m)} \in \mathbb{C}^{n\times n}$ is the discretized PDE Helmholtz operator \cite{Pratt_1998_GNF, Chen_2013_OFD} given by 
\begin{eqnarray} \label{A} 
\bold{A(m)} &= & \Delta + \omega^2 \bold{C}(\bold{m}) \text{diag}(\bold{m})\bold{B},
\label{helmholtz}
\end{eqnarray}
with $\omega$ the angular frequency and $\bold{\Delta}$ the discretized Laplace operator. 
$\bold{C}$ embeds boundary conditions and can be dependent or independent on $\bold{m}$ depending on the kinds of absorbing boundary conditions (radiation versus sponge)  \cite{Aghamiry_2019_IWR}. Also, $\bold{B}$ is used to spread the "mass" term $\omega^2\bold{C}(\bold{m}) \text{diag}(\bold{m})$ 
over all the  coefficients of the stencil to improve its accuracy following an anti-lumped mass strategy \cite{Marfurt_1984_AFF,Jo_1996_OPF,Hustedt_2004_MGS}.  \\
The second constraint in \eqref{main}, $\bold{Pu}= \bold{d}$, is the observation equation, in which $\bold{d} \in \mathbb{C}^{m\times 1}$
is the recorded seismic data and $\bold{P} \in \mathbb{R}^{m\times n}$ is a linear operator that samples the wavefield at the receiver positions.

\subsubsection{Reduced approach to solving \eqref{main}}
The reduced approach, which is more commonly use for sake of computational efficiency, strictly enforces the PDE constraint at each iteration by projection of the full space onto the parameter search space, leading to the following unconstrained optimization problem \cite{Pratt_1998_GNF,Plessix_2006_RAS,Virieux_2009_OFW}
\begin{equation} \label{pratt_FWI}
\underset{\bold{m}\in \mathcal{C}}{\min}~~  \Phi(\bold{m}) + \frac{\lambda}{2}||\bold{P}\bold{A}^{-1}(\bold{m})\bold{b}-\bold{d}||_2^2,
\end{equation}
where $\lambda>0$ is the penalty parameter. A number of methods have been proposed to solve the optimization problems of the form \eqref{pratt_FWI}, either for unregularized form $\Phi(\bold{m})=0$ \cite{Pratt_1998_GNF}, or the regularized form \cite{Asnaashari_2013_RSF,Guitton_2013_ISS,Esser_2018_TVR}. 
Although this reduced approach is more computationally tractable than the full-space approach, the highly-oscillating nature of  the inverse PDE operator $\bold{A}^{-1}$ makes the inverse problem highly nonlinear, and hence prone to convergence to a spurious local minima when the initial $\bold{m}$ is not accurate enough \cite{Symes_2008_MVA,Virieux_2009_OFW}.
The extended approach described below is an alternative way which is more immune to local minima. 
\subsubsection{WRI approach to solving \eqref{main}}
The extended approach, known as wavefield reconstruction inversion (WRI) \cite{VanLeeuwen_2013_MLM}, recasts the constrained optimization problem, equation \eqref{main}, as an unconstrained problem where both constraints are implemented with quadratic penalty functions.
\begin{equation} \label{leeuwen_FWI_closed}
\underset{\bold{u},\bold{m}\in \mathcal{C}}{\min}~~  \Phi(\bold{m})  + \frac{\lambda}{2}\|\bold{Pu}- \bold{d}\|_2^2 + \frac{\gamma}{2}\|\bold{A(m)u}-\bold{b}\|_2^2,
\end{equation}
where $\lambda,\gamma>0$ are the penalty parameters.
For the unregularized case where $\Phi(\bold{m})=0$ and also without the bounding constraint, \cite{VanLeeuwen_2013_MLM} solved this biconvex minimization problem with an alternating-direction algorithm, whereby the joint minimization over $\bold{u}$ and $\bold{m}$ is replaced by an alternating minimization over each variable separately.  
The main property of the penalty formulation given by equation \eqref{leeuwen_FWI_closed} is that
the PDE constraint in the original problem is replaced by a quadratic penalty term, which enlarges the search space and mitigates the inversion nonlinearity accordingly \cite{VanLeeuwen_2013_MLM}.
Its main drawback, however, is the difficulty related to the adaptive tuning of the penalty parameter, which is common to all penalty methods \cite{Nocedal_2006_NOO}.
\subsubsection{IR-WRI approach to solving \eqref{main}}
To overcome the above limitation, the iteratively-refined WRI (IR-WRI) implements the original constrained problem \eqref{main} 
with the augmented Lagrangian (AL) method \cite{Hestenes_1969_MUL,Nocedal_2006_NOO}. 
\begin{equation} \label{AL}
\begin{split}
\underset{\bold{u},\bold{m}\in \mathcal{C}}{\min} ~ \underset{\bold{v}_1,\bold{v}_2}{\max} ~ \Phi(\bold{m}) & + \frac{\lambda}{2}\|\bold{Pu}- \bold{d}\|_2^2 + \frac{\gamma}{2}\|\bold{A(m)u}-\bold{b}\|_2^2 \\
& + \bold{v}_1^T[\bold{Pu}- \bold{d}] +\bold{v}_2^T[\bold{A(m)u}-\bold{b}],
\end{split}
\end{equation}
where $\bold{v}_1$ and $\bold{v}_2$ are the dual variables (the Lagrangian multipliers).
The min max problem \eqref{AL} can also be written in a more compact form (the scaled form AL) as
\begin{equation} \label{ALc}
\begin{split}
\underset{\bold{u},\bold{m}\in \mathcal{C}}{\min} ~ \underset{\bold{v}_1,\bold{v}_2}{\max}  ~ \Phi(\bold{m}) & + \frac{\lambda}{2}\|\bold{Pu}- \bold{d}+\frac{1}{\lambda}\bold{v}_1\|_2^2-\frac{\lambda}{2}\|\bold{v}_1\|_2^2  \\
& + \frac{\gamma}{2}\|\bold{A(m)u}-\bold{b}+\frac{1}{\gamma}\bold{v}_2\|_2^2-\frac{\gamma}{2}\|\bold{v}_2\|_2^2,
\end{split}
\end{equation}
Applying a gradient ascent to \eqref{ALc} with respect to the duals, after a simple change of variables $\bold{d}^k=-\bold{v}^k_1/\lambda$ and $\bold{b}^k=-\bold{v}^k_2/\gamma$, gives the following iteration:
\begin{equation} \label{ALda}
\begin{split}
&\underset{\bold{u},\bold{m}\in \mathcal{C}}{\min}~~  \Phi(\bold{m})  + \frac{\lambda}{2}\|\bold{Pu}- \bold{d}- \bold{d}^k\|_2^2  \\
& \quad \quad \quad \quad\quad\quad\quad~ + \frac{\gamma}{2}\|\bold{A(m)u}-\bold{b}- \bold{b}^k\|_2^2,\\
&\bold{d}^{k+1} = \bold{d}^k + \bold{d} - \bold{Pu}, \\
&\bold{b}^{k+1} = \bold{b}^k + \bold{b} - \bold{A(m)u}, \\
\end{split}
\end{equation}
beginning with $\bold{d}^0=0$ and $\bold{b}^0=0$.
Capitalizing on the bilinearity of the wave equation in $\bold{m}$ and $\bold{u}$, ADMM \cite{Boyd_2011_DOS} is an efficient method to solve this kind of multivariate optimization problem. ADMM updates $\bold{m}$ and $\bold{u}$ separately through a Gauss-Seidel like iteration, i.e., fixing $\bold{m}$ and solving for $\bold{u}$ and vice versa. Accordingly, beginning with an initial guess $\bold{m}^0$, we end up with the following iteration to solve \eqref{ALda} \cite{Aghamiry_2019_IWR,Aghamiry_2018_MIA,Aghamiry_2019_IBCb}: 
\begin{subequations}
\label{main_ADMM}
\begin{align}
\bold{u}^{k+1} &=    \underset{\bold{u}}{\arg\min} ~~
\biggr\lVert
\begin{bmatrix}
\sqrt{\frac{\lambda}{\gamma}}\bold{P} \\
 \bold{A}(\bold{m}^{k})
\end{bmatrix} \bold{u} -
\begin{bmatrix}
\sqrt{\frac{\lambda}{\gamma}}(\bold{d}+\bold{d}^k)\\
\bold{b}+\bold{b}^k\\
\end{bmatrix} \biggr\rVert_2^2\label{main_ADMM_u}\\
\bold{m}^{k+1} &=  \underset{\bold{m} \in \mathcal{C}}{\arg\min} ~~ \Phi(\bold{m}) +  \frac{\gamma}{2}\|\bold{A(m)u}^{k+1}-\bold{b}- \bold{b}^k\|_2^2, \label{main_ADMM_m} \\
\bold{d}^{k+1} &= \bold{d}^k +\bold{d} - \bold{Pu}^{k+1}, \label{main_joint_d}\\
\bold{b}^{k+1} &= \bold{b}^k +\bold{b} - \bold{A}(\bold{m}^{k+1})\bold{u}^{k+1}, \label{main_joint_b}
\end{align}
\end{subequations} 
The subproblem \eqref{main_ADMM_u} associated with the wavefield reconstruction is quadratic and admits a closed-form solution. It relaxes the requirement to satisfy exactly the wave equation ($\bold{A}(\bold{m}^k)\bold{u}=\bold{b}$) for the benefit of an improved data fitting ($\bold{Pu=d}$). 
This is achieved by reconstructing the wavefields that best jointly fit the observations and satisfy the wave equation in a least-squares sense, while wavefields generated by the reduced approach exactly satisfy the wave equation, $\bold{u}_r=\bold{A}(\bold{m}^{k})^{-1}\bold{b}$, which makes the classical FWI highly non-convex \cite{Aghamiry_2019_IWR}.\\ 
Equation \eqref{ALda} shows that the duals are updated with the running sum of the data and source residuals in iterations and are used to update the right-hand sides in the penalty functions of the scaled AL. These error correction terms in the AL method are the key ingredients that allow for a constant penalty parameter to be used in iterations, while guaranteeing convergence to accurate minimizer \cite{Nocedal_2006_NOO}. \\
In the next section we focus on the solution of the model subproblem \eqref{main_ADMM_m} when compound regularizations are used as the regularization term.
\section{The model Subproblem}
This section presents different forms of the regularization functional including simple and compound functionals and the details of
our approach to solve the model subproblem \eqref{main_ADMM_m} with these functionals.

\subsection{Simple regularizers}
The two most widely used regularizers rely on the (squared) $\ell_2$ and $\ell_1$-norms.
The  squared $\ell_2$-norm, defined as
\begin{equation}
\|\bold{x}\|_2^2 = \sum_{i=1}^n \bold{x}_i^2,
\end{equation}
promotes smooth reconstruction, since the minimization of the squared value of components will penalize large components more severely than small ones.

In contrast, the $\ell_1$-norm, defined as
\begin{equation}
\|\bold{x}\|_1 = \sum_{i=1}^n |\bold{x}_i|,
\end{equation}
promotes sparse reconstruction (with many zero components), since the minimization of the absolute value of components will penalize   
small components more severely than the large counterparts.

The priors can be defined under a suitable transformation. 
For example, one may minimize the $\ell_1$- or $\ell_2$-norms of the first and/or second order differences of the model.
The first order forward differences in $x$- and $z$-direction are denotes by $\nabla_x f$ and $\nabla_z f$, with
\begin{equation}
\begin{cases}
(\nabla_x f)_{i,j} =f_{i,j} - f_{i,j-1},\\
(\nabla_z f)_{i,j} =f_{i,j} - f_{i-1,j},
\end{cases}
\end{equation}
with appropriate boundary conditions, where $i$ and $j$ run over the domain of the model parameters.
Accordingly, the discrete first order operator is defined as $\nabla=\begin{bmatrix}\nabla_x^T& \nabla_z^T \end{bmatrix}^T$  with
\begin{equation} \label{TV}
(|\nabla f|)_{i,j} = \sqrt{(\nabla_x f)_{i,j}^2 + (\nabla_z f)_{i,j}^2}.
\end{equation}
The $\ell_2$-norm of $|\nabla f|$ gives the first order Tikhonov regularization \cite{Tikhonov_1963_RIP}, which returns a flat regularized model (with a small gradient), while
its $\ell_1$-norm gives the total variation regularization \cite{Rudin_1992_NTV}, which returns a piecewise constant model (with a sparse gradient).

Analogously, the second order forward differences are denoted by $\nabla_{xx} f$ and $\nabla_{zz} f$, with
\begin{equation}
\begin{cases}
(\nabla_{xx} f)_{i,j} =f_{i,j-1} - 2f_{i,j}+f_{i,j+1}, \\
(\nabla_{zz} f)_{i,j} =f_{i-1,j} - 2f_{i,j}+f_{i+1,j},
\end{cases}
\end{equation}
with appropriate boundary conditions, where again $i$ and $j$ run over the domain of the model parameters.
Accordingly, the discrete second order operator is defined as $\nabla^2=\begin{bmatrix}\nabla_{xx}^T& \nabla_{zz}^T \end{bmatrix}^T$ with
\begin{equation} \label{TV2_1}
(|\nabla^2 f|)_{i,j} = \sqrt{(\nabla_{xx} f)_{i,j}^2 + (\nabla_{zz} f)_{i,j}^2}.
\end{equation}
The $\ell_2$-norm of $|\nabla^2 f|$ gives the second-order Tikhonov regularization, which returns a smooth regularized model (with a small Laplacian), while
its $\ell_1$-norm gives the second order TV regularization, which returns a piecewise linear model (with a sparse Laplacian).

Mixed second-order differences can also be constructed as $\nabla_{xz} f\equiv \nabla_z\nabla_x f$ with
\begin{equation}
(\nabla_{xz} f)_{i,j} =f_{i,j} -f_{i,j-1} - f_{i-1,j} +f_{i-1,j-1}.
\end{equation}
A discrete second-order operator, which includes mixed differences is defined as $\nabla^2=\begin{bmatrix}\nabla_{xx}^T& \sqrt{2}\nabla_{xz}^T& \nabla_{zz}^T \end{bmatrix}^T$ with 
\begin{equation} \label{TV2_2}
(|\nabla^2 f|)_{i,j} = \sqrt{(\nabla_{xx} f)_{i,j}^2+ 2(\nabla_{xy} f)_{i,j}^2 + (\nabla_{yy} f)_{i,j}^2},
\end{equation}
which equals the Frobenius norm of the Hessian matrix \cite{Lefkimmiatis_2012_HBN,Gholami_2013_FBS,Gholami_2019_3DD}.

Simple regularizations are effective for recovering models which can be characterize by a single prior. In the next section we propose compound regularizers which are more effective for recovering complicated models.

\subsection{Compound regularizers}
 Compound regularizers are constructed by combining two or more separate simple regularizers.
 This can be done by either a \textit{convex combination} (CC) or an \textit{infimal convolution} (IC).
\subsubsection{Convex combination of simple regularizers}
A CC of $r$ simple regularizer functionals $\Phi_1, ..., \Phi_{r}$ is
a compound regularizer functional of the form
\begin{equation}
\Phi_{\alpha}(\bold{x}) = \alpha_1\Phi_1(\bold{x})+ ...+ \alpha_r\Phi_{r}(\bold{x}),
\end{equation}
where weights $\alpha _{i}$ satisfy $\alpha _{i}\geq 0$ and 
\begin{equation}
\alpha_1+\alpha_2, ...,+ \alpha_r=1.
\end{equation}
Definitely, if all of the functions $\Phi_1, ..., \Phi_{r}$ are convex then $\Phi$ is so. In CC models, the regularized solution is forced to satisfy the individual priors simultaneously. 
As an example, a compound regularizer functional constructed by a CC of $\ell_1$- and squared $\ell_2$-norms ($\ell_1+\ell_2$), which known as elastic net \cite{Zou_2005_RAV,Gholami_2013_STF}, is 
\begin{equation} \label{CC}
\Phi_{\alpha}(\bold{x}) = (1-\alpha)\|\bold{x}\|_2^2 + \alpha\|\bold{x}\|_1,
\end{equation}
with $0\leq \alpha \leq 1$. 
The convexity of $\ell_1$- and $\ell_2$-norms implies that $\Phi_{\alpha}(\bold{x})$ in \eqref{CC} is convex.
One may also construct a compound regularizer functional by a CC of two $\ell_1$-norms which are applied in different domains, such as those spanned by two different wavelet transforms, or those spanned by a wavelet transform and the gradient operator \cite{Gholami_2010_ROL}.  

\subsubsection{Infimal convolution of simple regularizers}
In IC models, the solution is decomposed into simple components and then each component is regularized by an appropriate prior.
Accordingly, the IC of $r$ simple regularizer functionals $\Phi_1, ..., \Phi_{r}$ is a compound functional of the form
\begin{equation} \label{IC}
\Phi_{\alpha}(\bold{x}) = \min_{\bold{x}=\bold{x}_1+...+\bold{x}_r}\{\alpha_1\Phi_1(\bold{x}_1)+ ...+ \alpha_r\Phi_{r}(\bold{x}_r)\}.
\end{equation}
In the case of two functionals, $\Phi_{\alpha}$ in \eqref{IC} takes the form
\begin{equation} \label{IC2}
\Phi_{\alpha}(\bold{x}) = \min_{\bold{z}}\{(1-\alpha)\Phi_1(\bold{x}-\bold{z})+ \alpha\Phi_2(\bold{z})\},
\end{equation}
which is similar to the classical formula of convolution, and hence the term \textit{infimal convolution}.

The IC of $\ell_1$- and $\ell_2$-norms ($\ell_1\oplus\ell_2$) is 
\begin{equation} \label{IC12}
\Phi_{\alpha}(\bold{x})=\min_{\bold{z}}\{(1-\alpha)\|\bold{x}-\bold{z}\|_2^2 + \alpha\|\bold{z}\|_1\},
\end{equation}
which is a denoising problem and the solution of which is unique and is obtained simply by the well-known soft-threshold function \cite{Donoho_1995_DBS}: 
\begin{equation} \label{soft}
\bold{z} = \max \left(1 - \frac{\alpha}{2(1-\alpha)|\bold{x}|},0\right) \circ \bold{x},
\end{equation} 
Putting $\bold{z}$ from \eqref{soft} into \eqref{IC12} gives that
\begin{equation} \label{Huber}
\Phi_{\alpha}(\bold{x})= 
\begin{cases}
(1-\alpha)|\bold{x}|^2~~~\quad \text{if}~ |\bold{x}| \leq \frac{\alpha}{2(1-\alpha)} \\
\alpha|\bold{x}|-\frac{\alpha^2}{4(1-\alpha)} \quad\text{if}~ |\bold{x}| >\frac{\alpha}{2(1-\alpha)} 
\end{cases},
\end{equation}
which is nothing other than the Huber function \cite{Huber_1973_RRA}.
As seen, this function has a hybrid behavior: it has a quadratic behavior for small values of $|\bold{x}|$ and linear behavior for large values. 
The parameter $\frac{\alpha}{2(1-\alpha)}$ determines where the transition from quadratic to linear behavior takes place. 

\begin{figure}[!h]
\begin{center}
\includegraphics[scale=0.7]{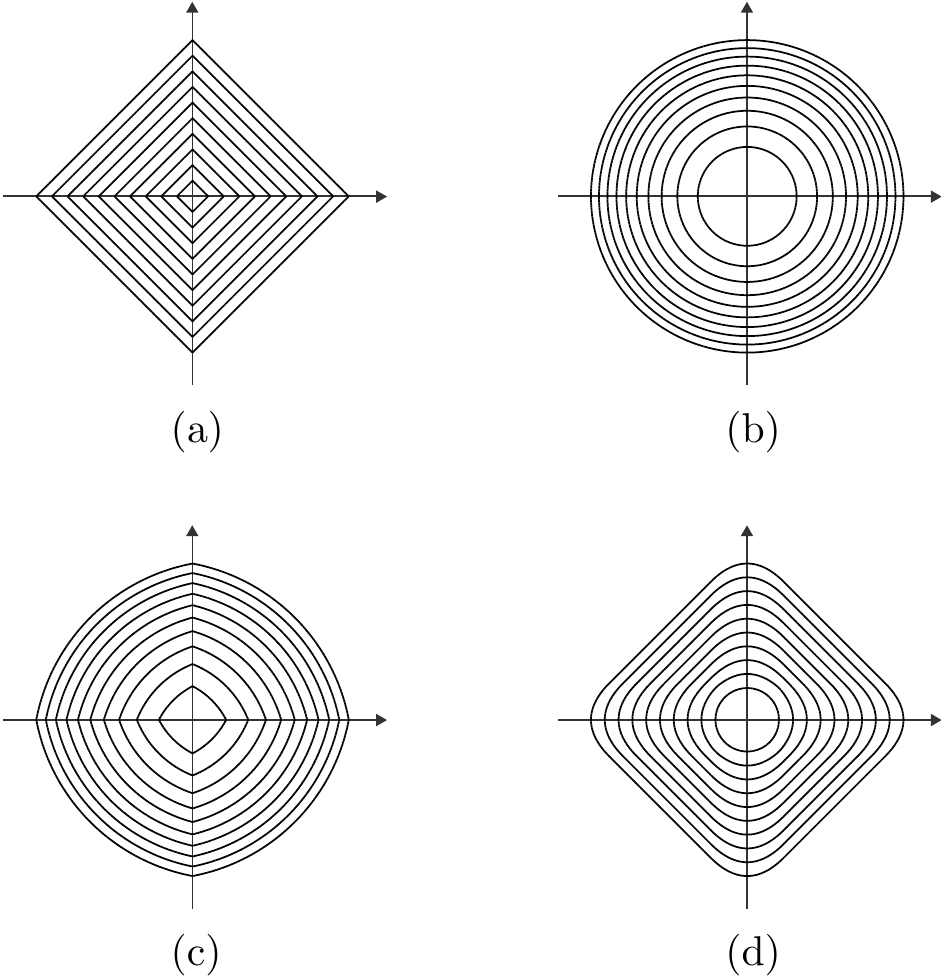}
\end{center}
\caption{Geometrical illustration of different regularizers. (a) the $\ell_1$-norm, (b) the $\ell_2$-norm, (c) the ($\ell_1+\ell_2$)-norm, and (d) the ($\ell_1\oplus\ell_2$)-norm.} \label{fig:GEOMs}
\end{figure}

Geometrical illustration of the $\ell_1$-norm, $\ell_2$-norm, ($\ell_1+\ell_2$)-norm, and ($\ell_1\oplus\ell_2$)-norm is shown in Fig. \ref{fig:GEOMs}. This figure shows that the $\ell_1$- and $\ell_2$-norms have a uniform behavior for all values, while the CC norm (the ($\ell_1+\ell_2$)-norm) has a hybrid behavior: it approaches the $\ell_1$-norm near zero, where it behaves as a linear function, but approaches the $\ell_2$-norm for large values,  where it behaves as a quadratic function. 
Unlike $\ell_1+\ell_2$, the IC function $\ell_1\oplus\ell_2$ approaches the $\ell_2$-norm near zero but is linear and approaches the $\ell_1$-norm for large values.

In this paper, we consider \eqref{IC2} in the following settings, though other configurations are possible:
\begin{equation} \label{TT}
\Phi^{\text{TT}}_{\alpha}(\bold{x}) = \min_{\bold{x}=\bold{x}_1+\bold{x}_2} (1-\alpha)\|\nabla^2\bold{x}_2\|_2^2+ \alpha\|\nabla\bold{x}_1\|_1,
\end{equation}
and
\begin{equation} \label{TGV}
\Phi^{\text{TGV}}_{\alpha}(\bold{x}) = \min_{\bold{x}=\bold{x}_1+\bold{x}_2} (1-\alpha)\|\nabla^2\bold{x}_2\|_1+ \alpha\|\nabla\bold{x}_1\|_1,
\end{equation}
where in both \eqref{TT} and \eqref{TGV}, the norms are applied on the absolute valued components of $\nabla^2\bold{x}_2$ (\eqref{TV2_1} and \eqref{TV2_2}) and $\nabla\bold{x}_1$ \eqref{TV}.
The compound regularizer $\Phi^{\text{TT}}_{\alpha}$ is a combination of the second order Tikhonov and TV (TT) regularizations  \cite{Gholami_2013_BCT} and $\Phi^{\text{TGV}}_{\alpha}$ is a combination of the first and second order TV regularizations, called total generalized variation (TGV) \cite{Bredies_2010_TGV,Setzer_2011_ICR}. 
The former is suitable for recovering piecewise-smooth models, while the latter is better suited for piecewise linear models.
Next section gives a solution procedure to solve \eqref{main_ADMM_m} with these regularizers.
%
%
%
%
%
%
\subsection{Solving the subproblem \eqref{main_ADMM_m}}
In this section we present how to solve the subproblem \eqref{main_ADMM_m} with TT regularization. The solution procedure for the TGV regularizer follows easily.  
From the definition of $\bold{A}$, in \eqref{A}, we get that
\begin{eqnarray} 
\bold{A}(\bold{m})\bold{u}  &=&  \Delta \bold{u}  + \bold{L} \bold{m},
\label{eqL}
\end{eqnarray}
where
\begin{equation}
 \bold{L}=\frac{\partial\bold{A(m)}}{\partial \bold{m}}\bold{u}=\omega^2 \bold{C} \text{diag}(\bold{B}\bold{u}),
\end{equation}
where we assume that $\bold{C}$ does not depend on $\bold{m}$ (this is the case for perfectly-matched absorbing boundary conditions  \cite{Aghamiry_2019_IWR}).
From the explicit decomposition $\bold{m}=\bold{m}_1+\bold{m}_2$ and \eqref{eqL}, the solution of the optimization problem \eqref{main_ADMM_m} can be expressed as 
\begin{equation}
 \underset{\substack{\bold{m}=\bold{m}_1+\bold{m}_2\\ \bold{m}\in  \mathcal{C}}}{\arg\min} ~~ \Phi^{\text{TT}}_{\alpha}(\bold{m}_1,\bold{m}_2) +  \frac{\gamma}{2}\|\bold{L} [\bold{m}_1+\bold{m}_2] - \bold{y}\|_2^2 , \label{MAIN}
\end{equation}
%
%
%
%
%
%
where $\bold{y}=\bold{b+b}^k-\Delta \bold{u}^{k+1}$.
Applying ADMM to \eqref{MAIN} leads to the following iteration
\begin{subequations}
\label{subsub}
\begin{align}
\begin{bmatrix}
\bold{m}_1^{k+1}\\
\bold{m}_2^{k+1}
\end{bmatrix}
=  \underset{\bold{m}_1,\bold{m}_2}{\arg\min} ~ 
C(\bold{m}_1,\bold{m}_2,{\bold{q}}^k,\tilde{\bold{q}}^k,\bold{p}^k,\tilde{\bold{p}}^k), \label{subsub_m}\\
 \bold{p}^{k+1} =  \underset{\bold{p}}{\arg\min} ~
   \alpha\|\bold{p}\|_1+\frac{\zeta}{2}  \| \bold{\nabla} \bold{m}_1^{k+1} -\bold{p}-\tilde{\bold{p}}^k\|_2^2,  \label{subsub_p}\\
 \bold{q}^{k+1} =  \underset{\bold{q} \in \mathcal{C}}{\arg\min} ~
\frac{\eta}{2}\|\bold{m}^{k+1}_1+\bold{m}^{k+1}_2-\bold{q}-\tilde{\bold{q}}^k\|_2^2, \label{subsub_q}
\end{align} 
\end{subequations}  
%
%
where
\begin{equation}
\begin{split}
& C(\bold{m}_1,\bold{m}_2,{\bold{q}}^k,\tilde{\bold{q}}^k,\bold{p}^k,\tilde{\bold{p}}^k)  =
 \frac{\gamma}{2}\|\bold{L}[\bold{m}_1+\bold{m}_2] -\bold{y} \|_2^2 \\ 
 & \qquad ~~~+ (1-\alpha)\|\nabla^2\bold{m}_2\|_2^2 +  \frac{\zeta}{2}\|\bold{\nabla}{\bold{m}_1}-\bold{p}^k-\tilde{\bold{p}}^k\|_2^2 \\
 & \hspace{3.1cm}~~~+\frac{\eta}{2}\|\bold{m}_1+\bold{m}_2-\bold{q}^k-\tilde{\bold{q}}^k\|_2^2,
\end{split}
\end{equation} 
and $\bold{p} = \nabla \bold{m}_1$ and $\bold{q}$ are auxiliary primal variables that are introduced to decouple the 
$\ell_1$ and the $\ell_2$ minimization problems and solve the former ones with proximal algorithms following the split Bregman method \cite{Goldstein_2009_SBM}.
The dual variables  $\tilde{\bold{g}}$, $\tilde{\bold{m}}$ are updated through a gradient ascent step according to the method of multipliers \cite{Nocedal_2006_NOO}
\begin{subequations} 
\begin{eqnarray}
\tilde{\bold{p}}^{k+1} &=& \tilde{\bold{p}}^k +\bold{p}^{k+1} -\bold{\nabla}\bold{m}_1^{k+1},\\
\tilde{\bold{q}}^{k+1} &=& \tilde{\bold{q}}^k +\bold{q}^{k+1} - [\bold{m}_1^{k+1}+\bold{m}_2^{k+1}],  
\end{eqnarray} 
\end{subequations}
%
%
%
%
We now discuss how to solve the subsubproblems given in \eqref{subsub}.

\subsubsection{The subsubproblem \eqref{subsub_m}}
A solution of subsubproblem \eqref{subsub_m} occurs at the point where the derivatives of the objective function $C$ with respect to $\bold{m}_1$ and $\bold{m}_2$ vanish simultaneously. 
Accordingly, we end up with  the following linear system of equations:

\begin{equation} \label{jointm}
\begin{bmatrix}
     \bold{G}_{11}  & \bold{G}_{12}  \\
   \bold{G}_{21}   &  \bold{G}_{22}   
\end{bmatrix}
\begin{bmatrix}
   \bold{m}_{1}  \\
   \bold{m}_{2}
\end{bmatrix}
=
\begin{bmatrix}
   \bold{h}_{1}\\
   \bold{h}_{2}  
\end{bmatrix},
\end{equation}
with
\begin{equation*}
\begin{cases}
\bold{G}_{11}=  \gamma \bold{L}^T\bold{L}+\zeta {\nabla}^T \nabla+\eta \bold{I},\\
\bold{G}_{12}= \bold{G}_{21}=\gamma   \bold{L}^T\bold{L}+\eta \bold{I},\\
\bold{G}_{22}= \gamma  \bold{L}^T\bold{L}+(1-\alpha) (\nabla^2)^{T} \nabla^2+\eta \bold{I},\\
\end{cases}
\end{equation*}
and
\begin{equation*}
\begin{cases}
\bold{h}_{1}=  \gamma  \bold{L}^T\bold{y}+\zeta \nabla^T[\bold{p}^k+\tilde{\bold{p}}^k]+\eta [\bold{q}^k+\tilde{\bold{q}}^k],\\
\bold{h}_{2}= \gamma  \bold{L}^T\bold{y}+\eta [\bold{q}^k+\tilde{\bold{q}}^k],
\end{cases}
\end{equation*}
where  $\bold{I}$ is the identity matrix. 

\cite{Aghamiry_2018_HTT} broke down the $2n\times2n$ problem \eqref{jointm} into two smaller $n \times n$ systems and updates $\bold{m}_1$ and $\bold{m}_2$ in alternating mode at the expense of convergence speed \cite[their eqs. 10 and 11]{Aghamiry_2018_HTT}.
Instead, we solve here the original system exactly using a variable projection scheme, thus leading to faster convergence and more accurate results. From the first equation of \eqref{jointm}, we find that 
\begin{equation} \label{vpl1}
\bold{m}_2=\bold{G}_{12}^{-1}[\bold{h}_1-\bold{G}_{11}\bold{m}_1]
\end{equation}
and plugging this into the second equation of \eqref{jointm} we get the following:
\begin{equation} \label{jm1m2}
 \bold{m}_1= [\bold{G}_{11}-\bold{G}_{22} \bold{G}_{12}^{-1} \bold{G}_{11}]^{-1}[\bold{h}_2-\bold{G}_{22}\bold{G}_{12}^{-1}\bold{h}_{1}].
\end{equation}
Interestingly, $\bold{L}$ is diagonal, implying that $\bold{G}_{12}$ is also diagonal.
Thus we only need to solve an $n\times n$ system to estimate $\bold{m}_1$, from which $\bold{m}_2$ easily follows.
\subsubsection{The subsubproblem \eqref{subsub_p}}
The sub-problem for $\bold{p}$, equation \eqref{subsub_p}, is a denoising problem and is straightforward to solve.
Note that $\bold{p}$ has two components associated with the gradient in each direction:
\begin{equation}
\bold{p} = \begin{bmatrix}
\bold{p}_x\\
\bold{p}_z
\end{bmatrix}.
\end{equation}
Equation \eqref{subsub_p} is solved with a generalized proximity operator \cite{Combettes_2011_PRO} leading to
\begin{equation}
\bold{p}^{k+1}=\text{prox}_{\zeta/\alpha}(\bold{z})=
\begin{bmatrix}
\xi\circ \bold{z}_x\\
\xi\circ \bold{z}_z
\end{bmatrix},
\end{equation}
where
\begin{equation}
\bold{z} = \nabla\bold{m}_1^{k+1}-\tilde{\bold{p}}^k=
\begin{bmatrix}
\bold{z}_x\\
\bold{z}_z
\end{bmatrix}, 
\end{equation}
and
\begin{equation} \label{prox0}
\xi = \max(1 - \frac{\zeta}{\alpha \sqrt{\bold{z}_x^2+\bold{z}_z^2}},0).
\end{equation}

\subsubsection{The subsubproblem \eqref{subsub_q}}
The optimization problem \eqref{subsub_q} has also a entrywise solution given by
 \begin{equation} \label{projp0}
 \bold{q}^{k+1} = \text{proj}_{\mathcal{C}} (\bold{m}^{k+1}_1+\bold{m}^{k+1}_2-\tilde{\bold{q}}^k),
 \end{equation}
where the projection operator projects its argument onto the desired box $[\bold{m}_l,\bold{m}_u]$
according to $\text{proj}_{\mathcal{C}} (\bullet) = \min(\max(\bullet,\bold{m}_{l}), \bold{m}_{u})$.

It should be noted that the total algorithm consists of two levels of iterations: an outer iteration given in \eqref{main_ADMM} and an inner iteration given in \eqref{subsub} corresponding to the model subproblem \eqref{main_ADMM_m}.
Numerical results, however, show that only one inner iteration suffices for convergence of the algorithm, hence significantly reducing the total computational cost.

\section{Numerical examples}
We assess the performance of our algorithm against 1D and 2D mono-parameter synthetic examples. 
In Table \ref{Tab:reg} we give different regularization functions which are applied for stabilizing the FWI solution.
 We start with zero-offset Vertical-Seismic-Profiling (VSP) examples (1D IR-WRI) where the targeted wave speed profiles are selected from well-documented 2D benchmark subsurface velocity models in exploration seismic. To tackle more realistic applications, we proceed with a target of the 2D challenging 2004 BP salt model \cite{Billette_2004_BPB} when a crude initial model and realistic frequencies are used as starting points.
\subsection{Performance comparison using 1D test on benchmark models}
\begin{figure}[!h]
\begin{center}
\includegraphics[scale=0.7]{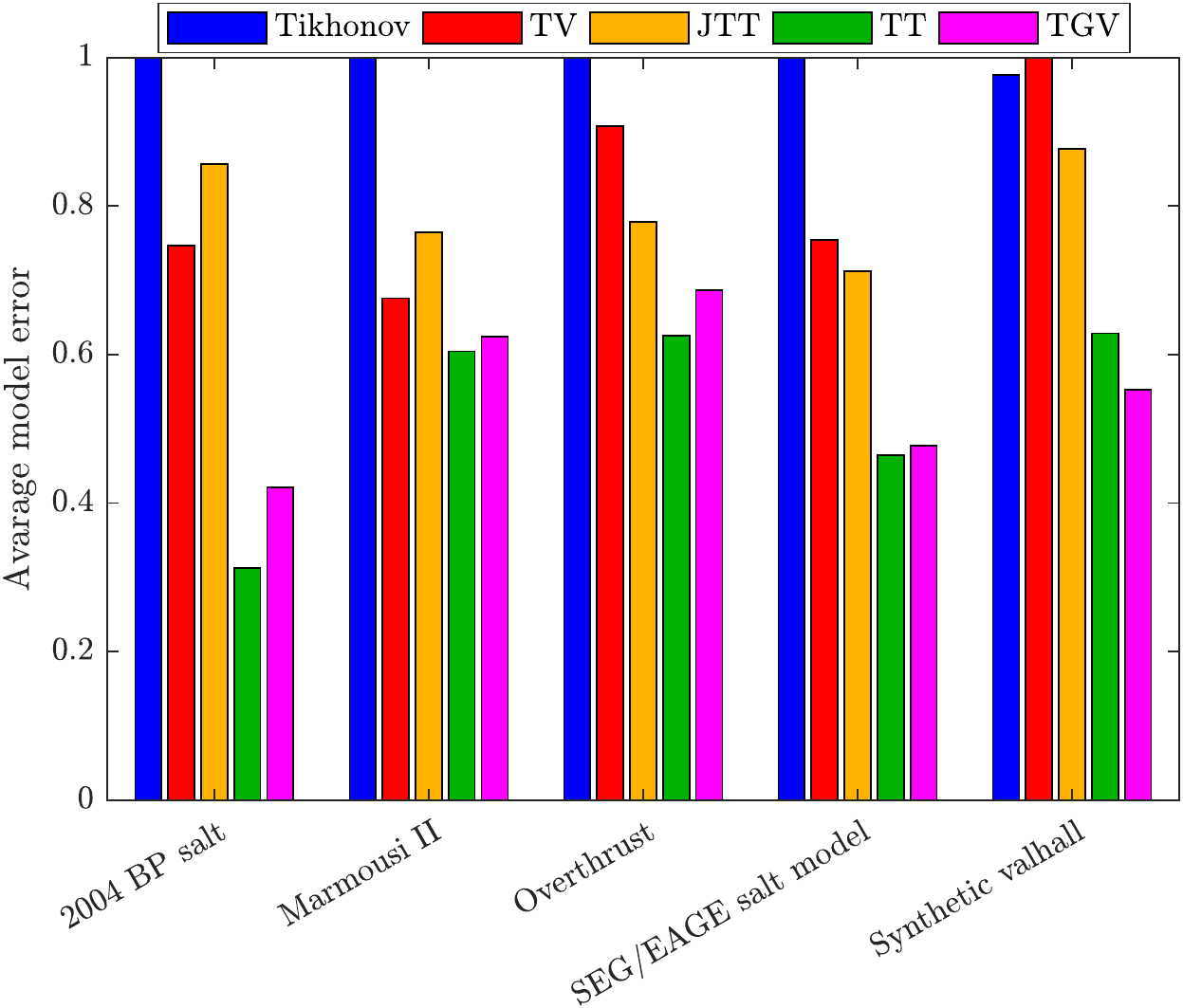}
\end{center}
\caption{Zero offset VSP test. Average model error in estimated 1D profiles for different velocity models and different regularization functions.} \label{fig:MSElog}
\end{figure}

\begin{figure} 
\centering
   \begin{subfigure}[b]{0.48\textwidth}
   \includegraphics[width=1\textwidth]{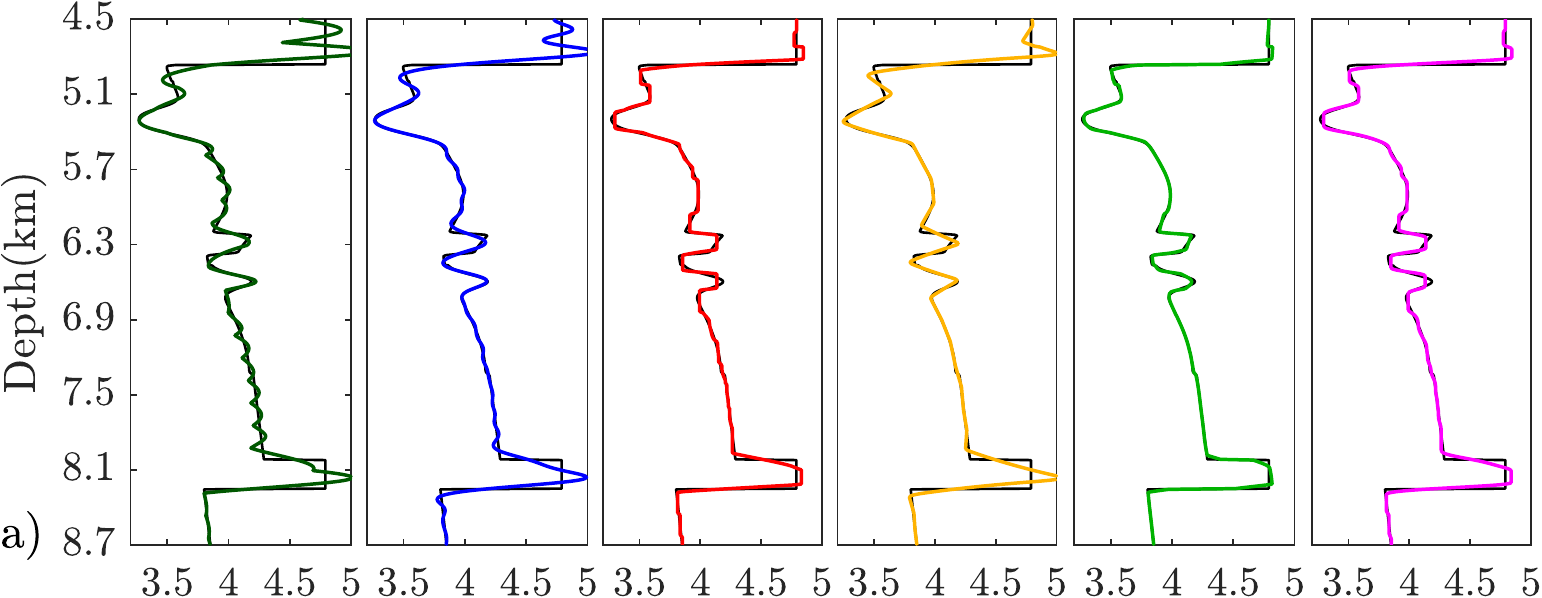} 
\end{subfigure}
\begin{subfigure}[b]{0.48\textwidth}
   \includegraphics[width=\textwidth]{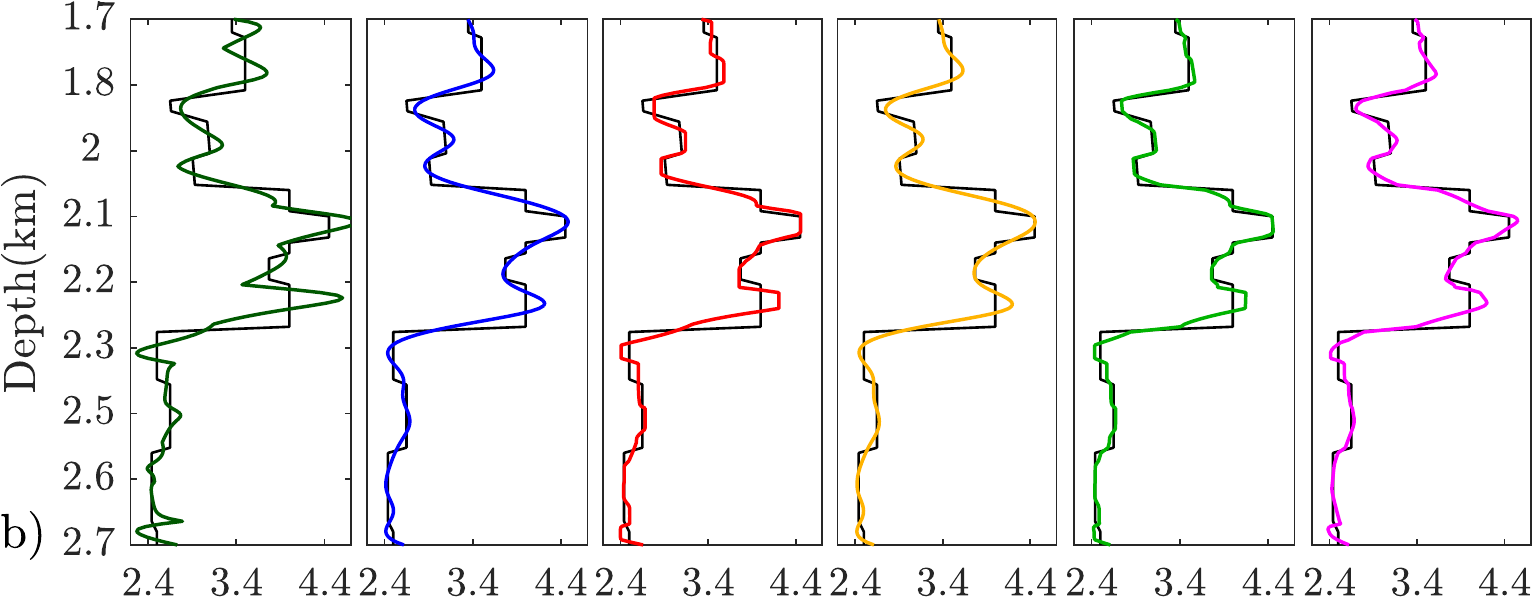}
\end{subfigure}
\begin{subfigure}[b]{0.48\textwidth}
   \includegraphics[width=\textwidth]{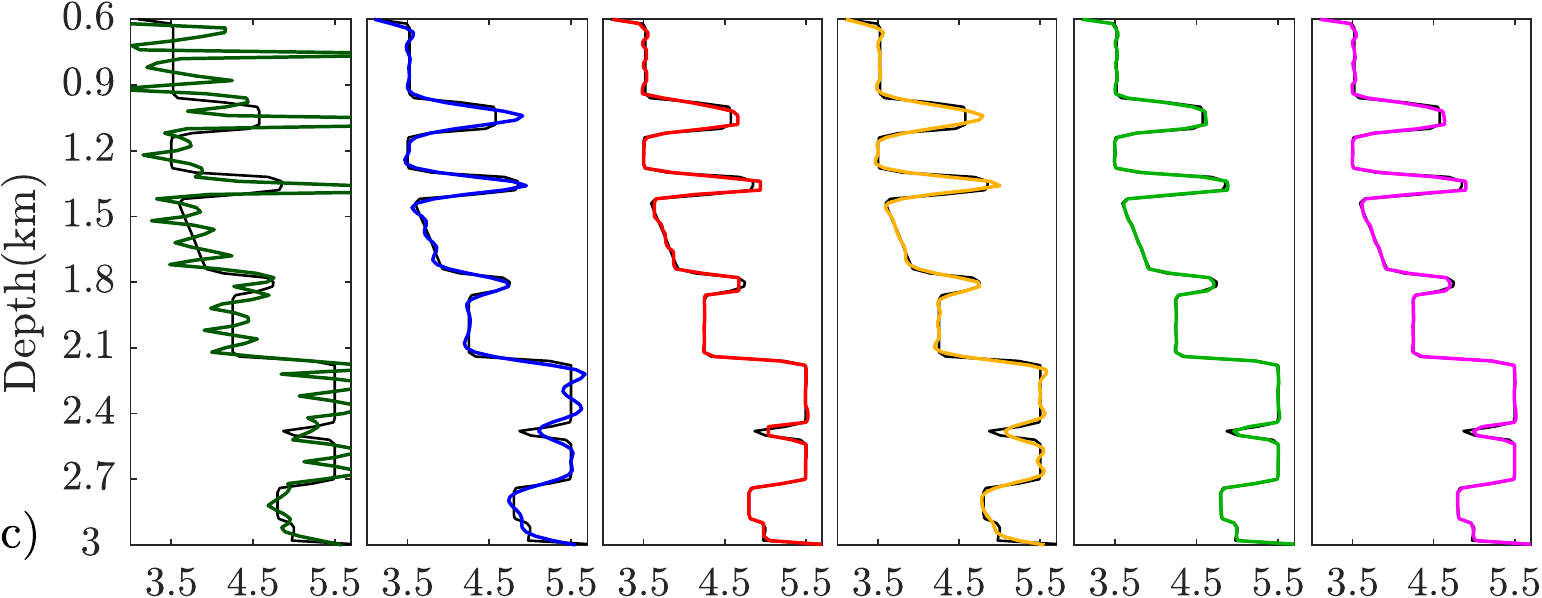}
\end{subfigure}
\begin{subfigure}[b]{0.48\textwidth}
   \includegraphics[width=\textwidth]{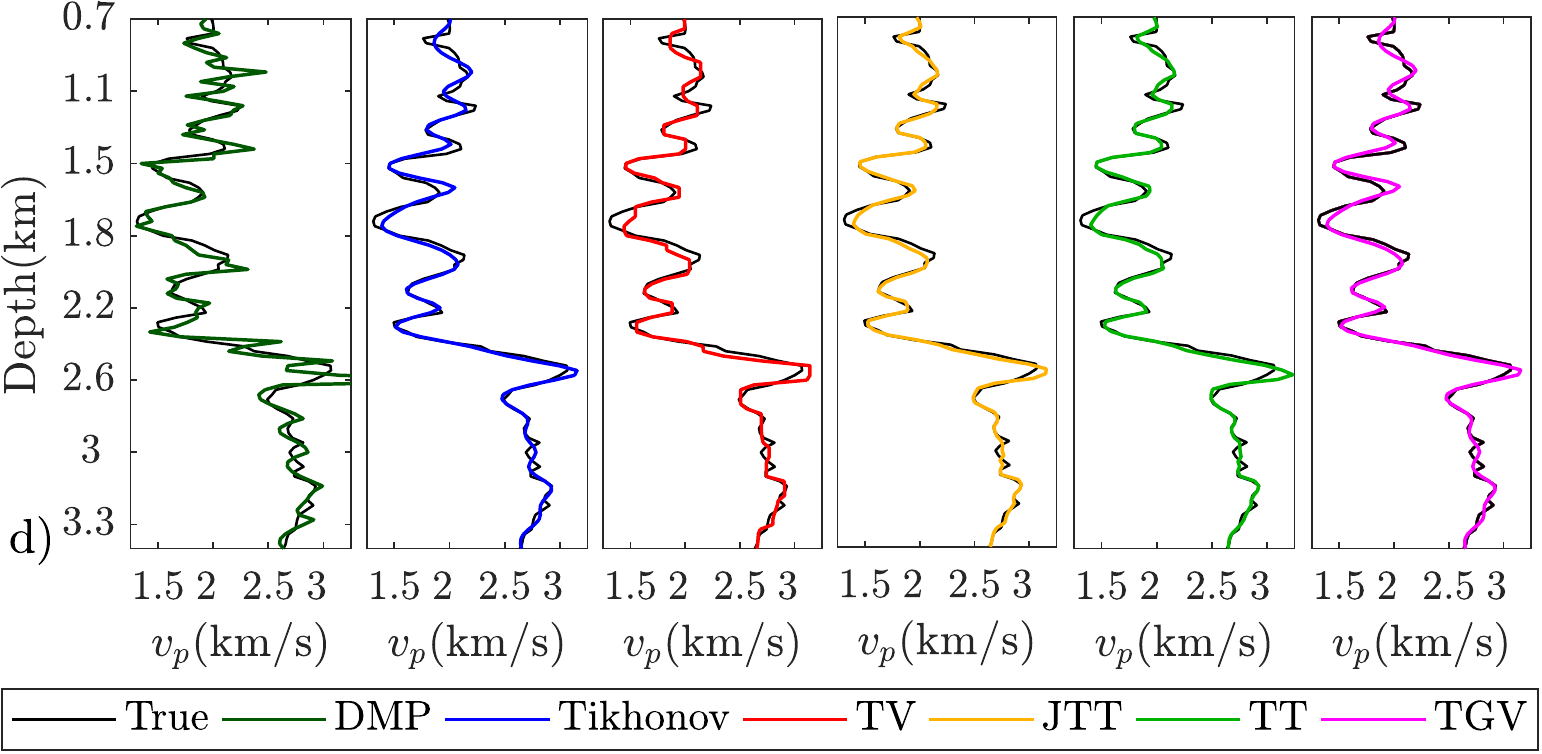}
\end{subfigure}
\caption{Zero offset VSP test. Some part of (a) 2004 Bp salt, (b) Marmousi II, (c) Overthrust, and (d) synthetic Valhall  models estimated  with different regularizations.}
\label{fig:part_log}
\end{figure}

First, we assess the performance of our regularized IR-WRI against 1D mono-parameter synthetic examples when the true models are 100 vertical profiles selected from the 2004 BP salt \cite{Billette_2004_BPB}, Marmousi II \cite{Martin_2006_M2E}, SEG/EAGE overthrust \cite{Aminzadeh_1997_DSO}, SEG/EAGE salt \cite{Aminzadeh_1997_DSO} and synthetic Valhall \cite{Prieux_2011_FAI} benchmark velocity models (we extracted 20 profiles from each benchmark model). For all of the experiments, a single source is used at the surface and a single frequency is considered for inversion. The model dimension, the inverted frequency and the receiver spacing are outlined for each model in Table \ref{Tab:log}. We perform forward modeling with a 3-point finite-difference stencil and PML absorbing boundary conditions at the two ends of the model. The starting model for IR-WRI is a homogeneous velocity model in which the velocity is the mean value of each profile.
We set the penalty parameters  according to the guideline proposed by \cite{Aghamiry_2019_IBCb} for bound-constrained TV regularized IR-WRI. Moreover, for a fair comparison of the compound regularizers (JTT, TT and TGV), we select for each of them the optimum value of $\alpha$ among a range of preset values that minimizes the error in the models estimated by the IR-WRI. Also, we set the parameter bounds $\bold{m}_l$ and $\bold{m}_u$ equal to 50\% and 150\% of the minimum and maximum velocities of the true model, respectively.      
\begin{table}[]
\caption{Different regularization function.}
\label{Tab:reg}
\begin{center}
\begin{tabular}{|c|c|}
\hline
Nick name                 &  Expression of $\Phi(\bold{m})$ \rule{0pt}{10pt} \\ [3pt]\hline
DMP          & $\|\bold{m}\|_2^2$     \rule{0pt}{10pt} \\ [3pt]\hline
Tikhonov   &    $\|\nabla^2\bold{m}\|_2^2$        \rule{0pt}{10pt} \\ [3pt]\hline
TV             & $\|\nabla\bold{m}\|_1$        \rule{0pt}{10pt} \\ [3pt]\hline
JTT       &       $(1-\alpha)\|\nabla^2\bold{m}\|_2^2+ \alpha\|\nabla\bold{m}\|_1$   \rule{0pt}{10pt}\\ [3pt]\hline
TT        &     $\underset{\bold{m}=\bold{m}_1+\bold{m}_2}{\min}\{(1-\alpha)\|\nabla^2\bold{m}_2\|_2^2+ \alpha\|\nabla\bold{m}_1\|_1\}$     \rule{0pt}{10pt} \\ \hline
TGV                &      $\underset{\bold{m}=\bold{m}_1+\bold{m}_2}{\min}\{(1-\alpha)\|\nabla^2\bold{m}_2\|_1+ \alpha\|\nabla\bold{m}_1\|_1\}$   \rule{0pt}{10pt} \\ [3pt]\hline
\end{tabular}
\end{center}
\end{table}
%
%
%
\begin{table}[]
\caption{Experimental setup of 1D model tests}
\label{Tab:log}
\begin{tabular}{c|c|c|c|c|}
\cline{2-5}
                                          & \begin{tabular}[c]{@{}c@{}}Length \\ (km)\end{tabular} & \begin{tabular}[c]{@{}c@{}}Inverted \\  frequency\\ (Hz)\end{tabular} & \begin{tabular}[c]{@{}c@{}}Grid \\ interval \\ (m)\end{tabular} & \begin{tabular}[c]{@{}c@{}}Receiver\\ interval\\  (m)\end{tabular} \\ \hline
\multicolumn{1}{|c|}{2004 BP salt}        & 11.46                                                  & 5                                                                     & 6                                                               & 180                                                                \\ \hline
\multicolumn{1}{|c|}{Marmousi II}         & 3.75                                                   & 12                                                                    & 5                                                               & 85                                                                 \\ \hline
\multicolumn{1}{|c|}{Overthrust}          & 4.6                                                    & 12                                                                    & 20                                                              & 120                                                                \\ \hline
\multicolumn{1}{|c|}{SEG/EAGE salt model} & 4.2                                                    & 10                                                                    & 20                                                              & 120                                                                \\ \hline
\multicolumn{1}{|c|}{Synthetic valhall}   & 5.22                                                   & 5                                                                     & 25                                                              & 175                                                                \\ \hline
\end{tabular}
\end{table}
The monochromatic inversion is performed with noiseless data using a maximum number of iteration equal to 100 as stopping criterion. The average mean error of the estimated velocity profiles for the five benchmark models and the different regularizations are plotted in Fig. \ref{fig:MSElog}. The errors in each model for different regularizations are normalized to 1 for sake of clarity (the error of DMP regularizer is not shown because of its worse performance). Fig. \ref{fig:MSElog} clearly shows that the compound regularizations based upon infimal convolution (TT and TGV) always behave better than the CC regularization and the single regularization functionals (TV and Tikhonov). To emphasize the effects of the different regularization functions, we plot some close-up of the reconstructed profiles in Fig. \ref{fig:part_log} as well as the profiles reconstructed the DMP regularization in the first column. These results show that TT provides the most accurate reconstruction for the 2004 BP salt (Fig. \ref{fig:part_log}a) and Overthrust (Fig. \ref{fig:part_log}c) models. This is consistent with the fact that the velocity trends of these two models match well the piecewise smooth prior. In contrast, TGV behaves slightly better than TT for the Valhall model, whose velocity trend is the closest one to the piecewise linear prior (Fig. \ref{fig:part_log}d). For Marmousi II (Fig. \ref{fig:part_log}b), TT and TGV give similar results. 
%
%
%
%
%
%
\subsection{2004 BP salt model}
\begin{figure}[!h]
\begin{center}
\includegraphics[scale=0.55]{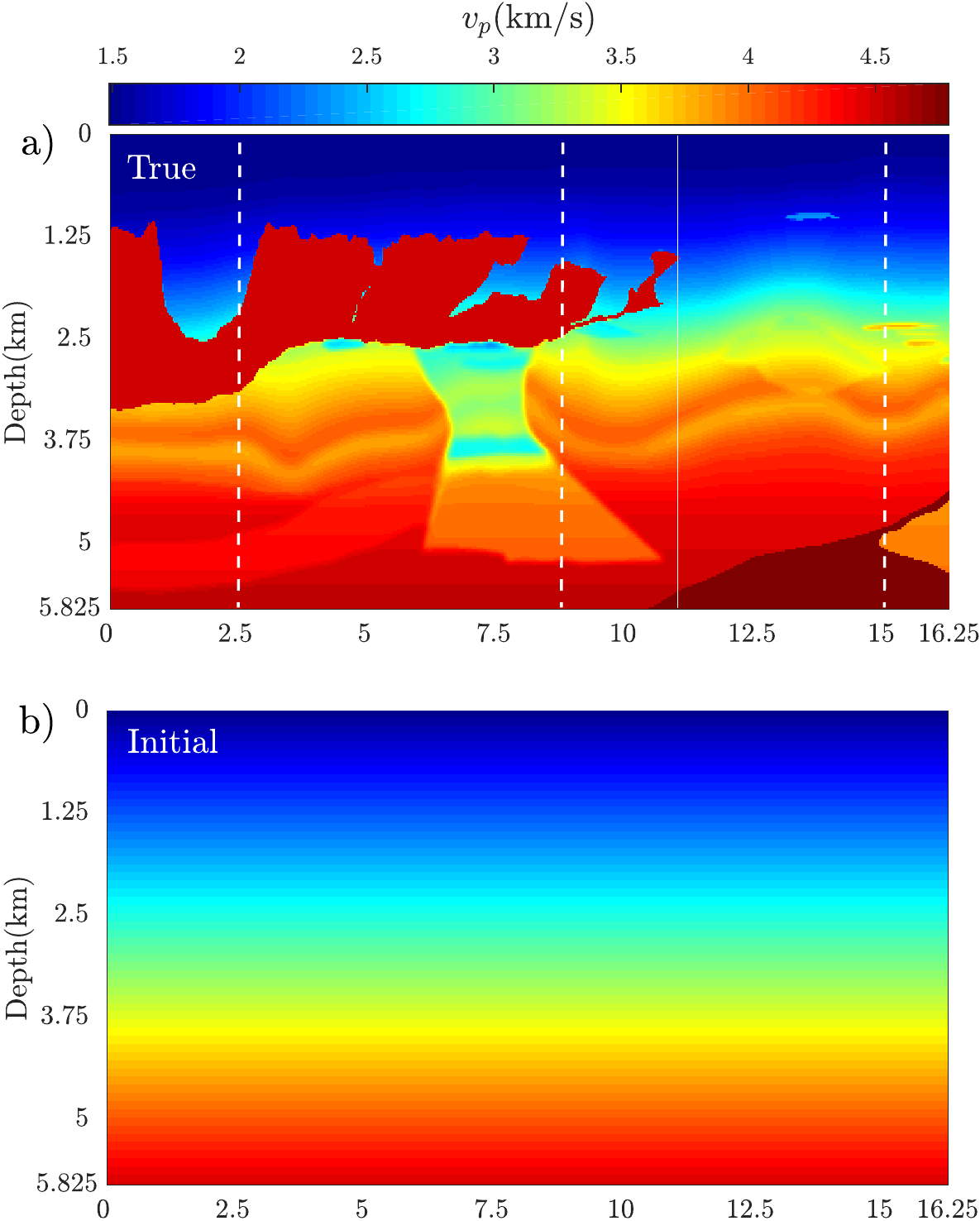}
\end{center}
\caption{2004 BP salt case study. (a) True model. The vertical dashed lines indicate the location of vertical logs of Fig.~\ref{fig:bp_start_log} and \ref{fig:bp_log}. (b) The velocity-gradient initial model.}
 \label{fig:bp_true}
\end{figure}

\begin{figure}[!h]
\begin{center}
\includegraphics[scale=0.38]{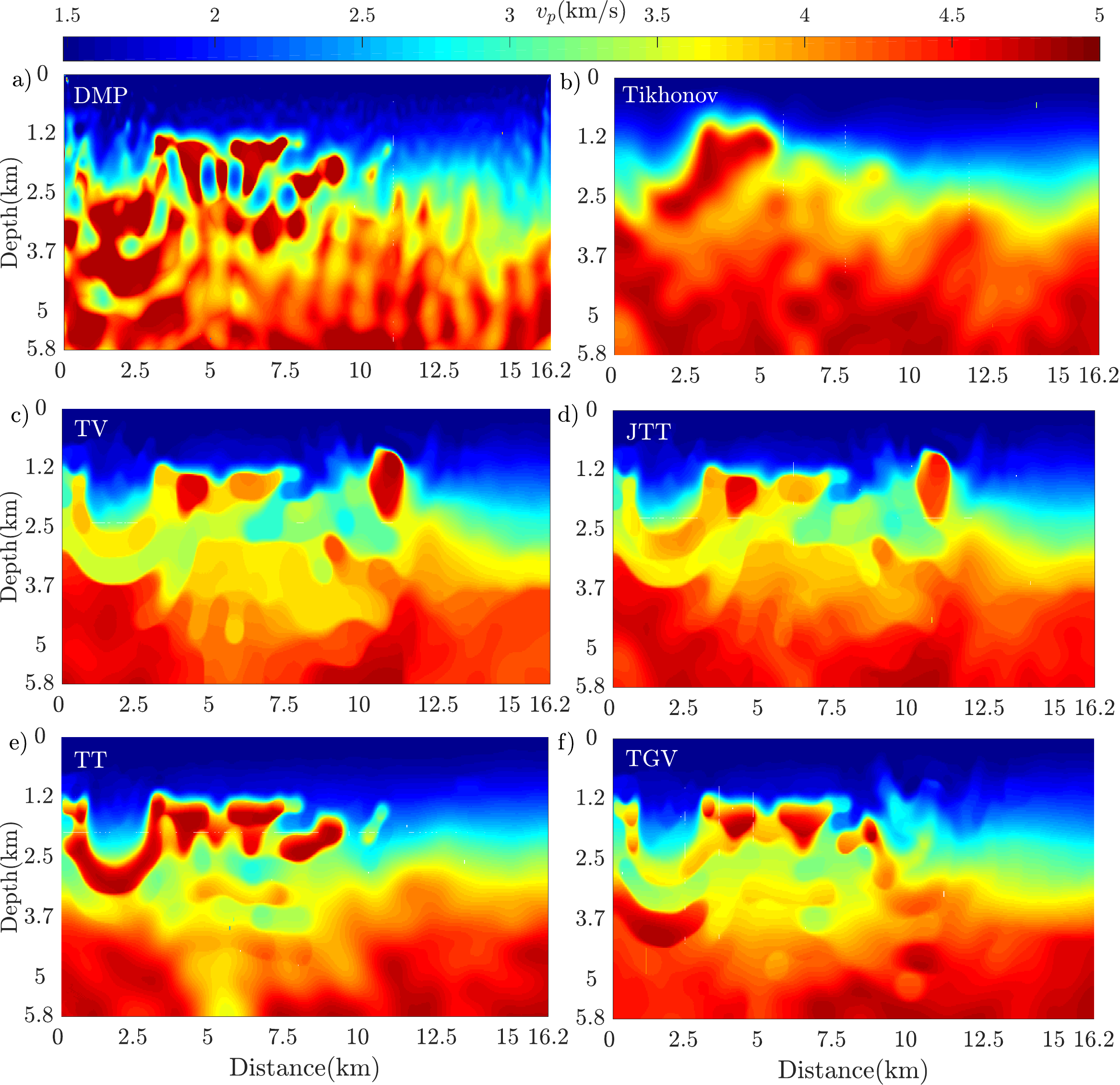}
\end{center}
\caption{2004 BP salt case study. \{3,3.5\}~Hz inversion when the velocity-gradient model (Fig. \ref{fig:bp_true}b) is used as initial model. (a-f) bound constrained IR-WRI with (a) DMP, (b) Tikhonov, (c) TV, (d) JTT, (e) TT and (f) TGV regularization.}
 \label{fig:bp_start}
\end{figure}

\begin{figure}[!h]
\begin{center}
\includegraphics[scale=0.45]{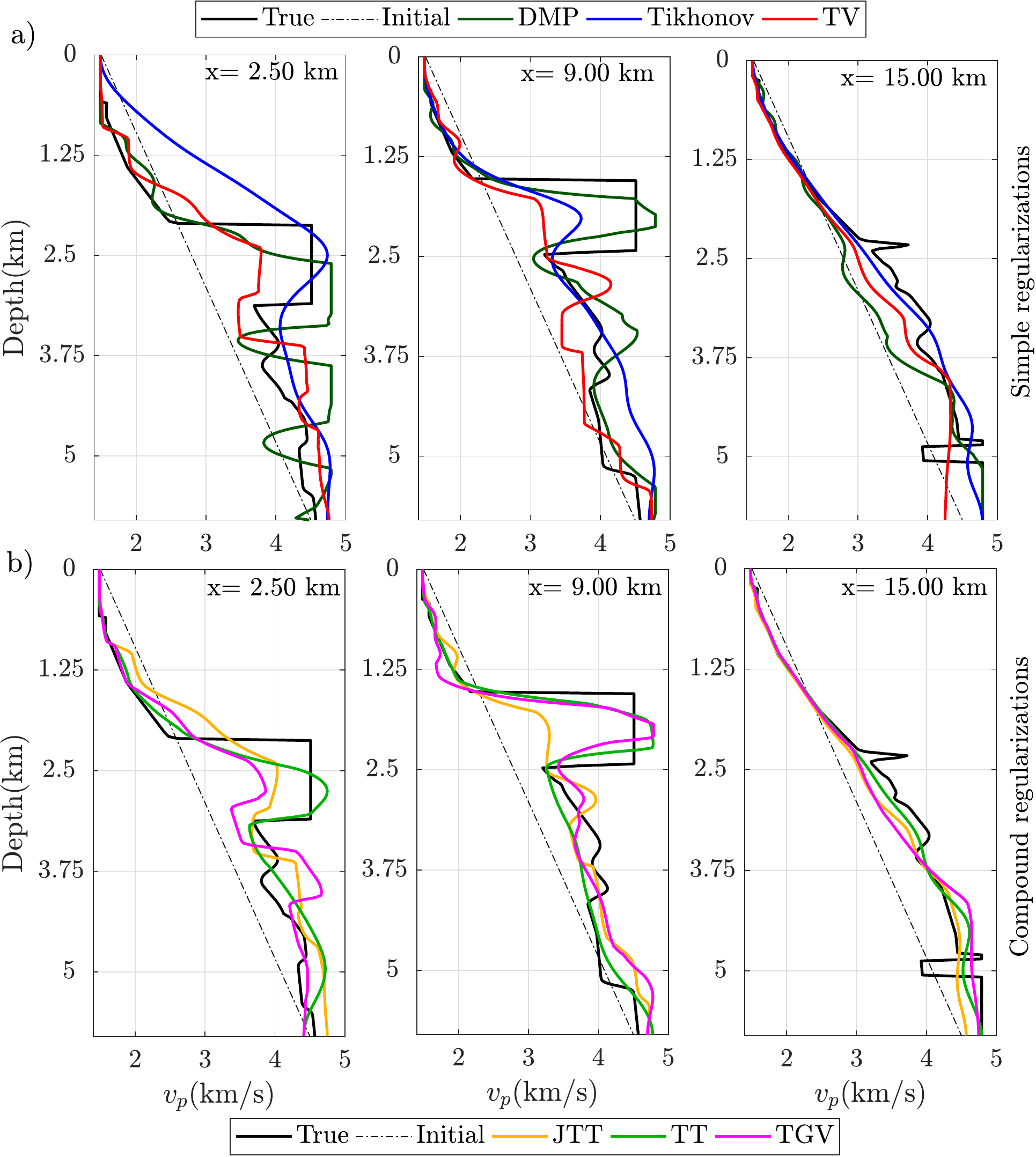}
\end{center}
\caption{2004 BP salt case study. \{3,3.5\}~Hz inversion with the velocity-gradient initial model. 
Direct comparison (along the logs shown in Fig.~\ref{fig:bp_true}a)  between the true velocity model (black), the initial model (dashed line) and the estimated models with (a) simple  regularizations (DMP in olive-green, Tikhonov in blue and TV in red) and
(b) compound regularizations (JTT in orange, TT in green and TGV in pink).}
 \label{fig:bp_start_log}
\end{figure}

\begin{figure} 
\centering
   \begin{subfigure}[b]{0.48\textwidth}
   \includegraphics[width=1\textwidth]{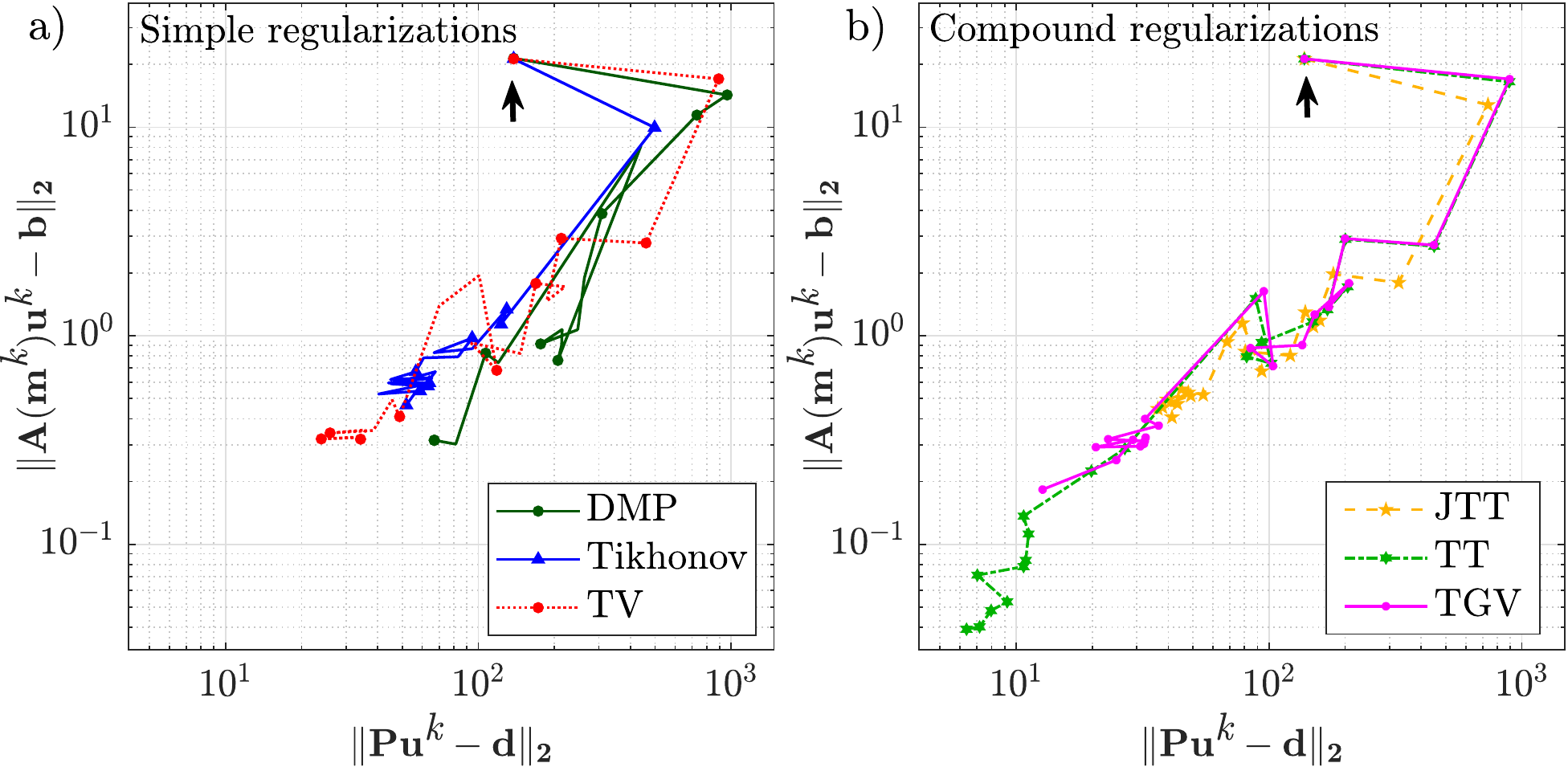} 
\end{subfigure}
\begin{subfigure}[b]{0.48\textwidth}
   \includegraphics[width=\textwidth]{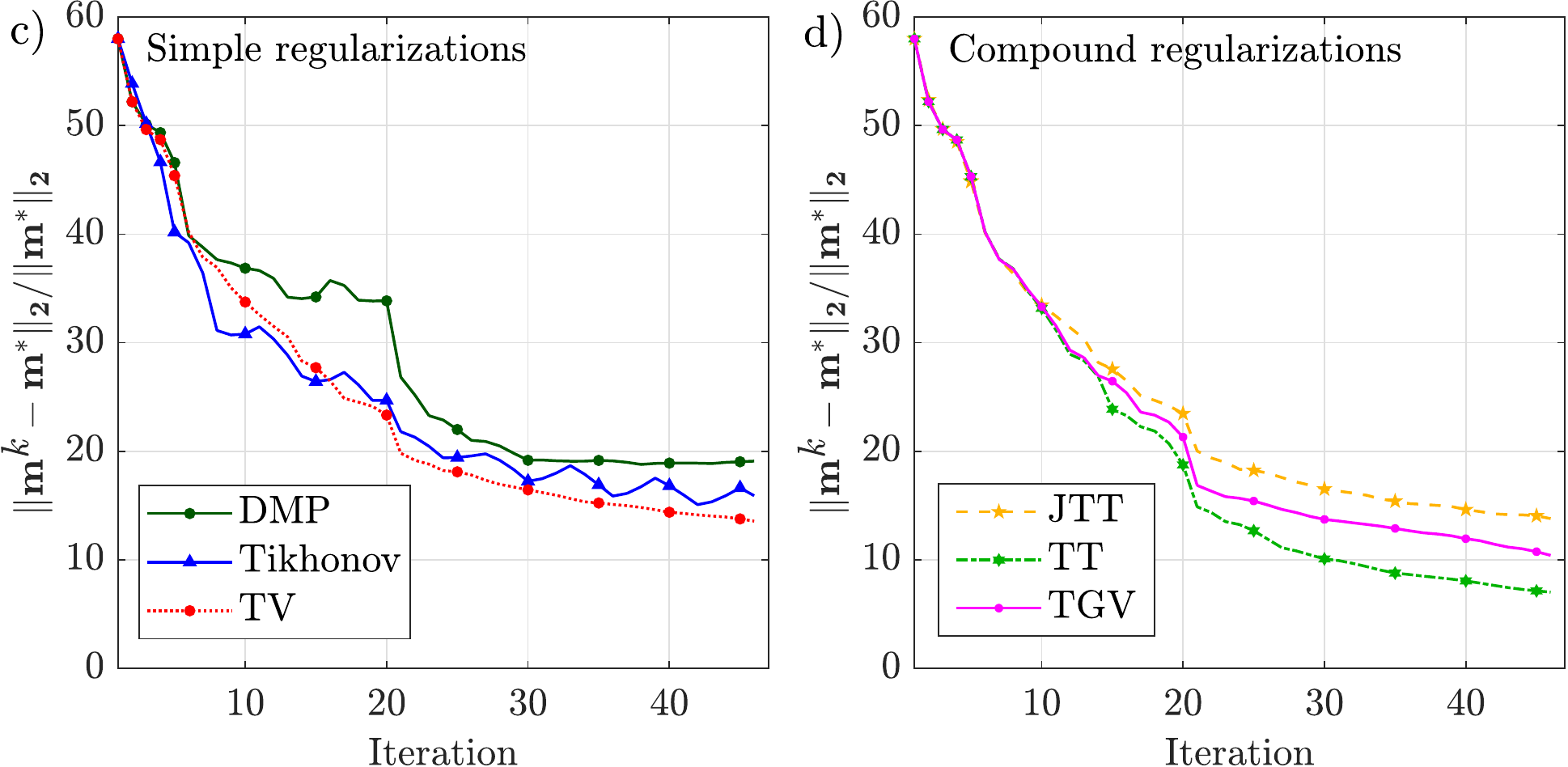}
\end{subfigure}
\caption{2004 BP salt case study. \{3,3.5\}~Hz inversion with the velocity-gradient initial model. 
(a-b) convergence history of the algorithm in the ($\|\bold{Pu}^k-\bold{d}\|_2 - \|\bold{A(m}^k)\bold{u}^k-\bold{b}\|_2$) plane for (a) simple  regularizations and (b) compound regularizations. The black arrow points the starting point. 
(c-d) evaluation of $\|\bold{m}^k-\bold{m}^*\|_2/\|\bold{m}^*\|_2$ during the iteration where $\bold{m}^*$ is the true model. The panels (a) and (b) as well as (c) and (d) are plotted with the same horizontal and vertical scale. }
\label{fig:bp_start_res}
\end{figure}

\begin{figure}[!h]
\begin{center}
\includegraphics[scale=0.38]{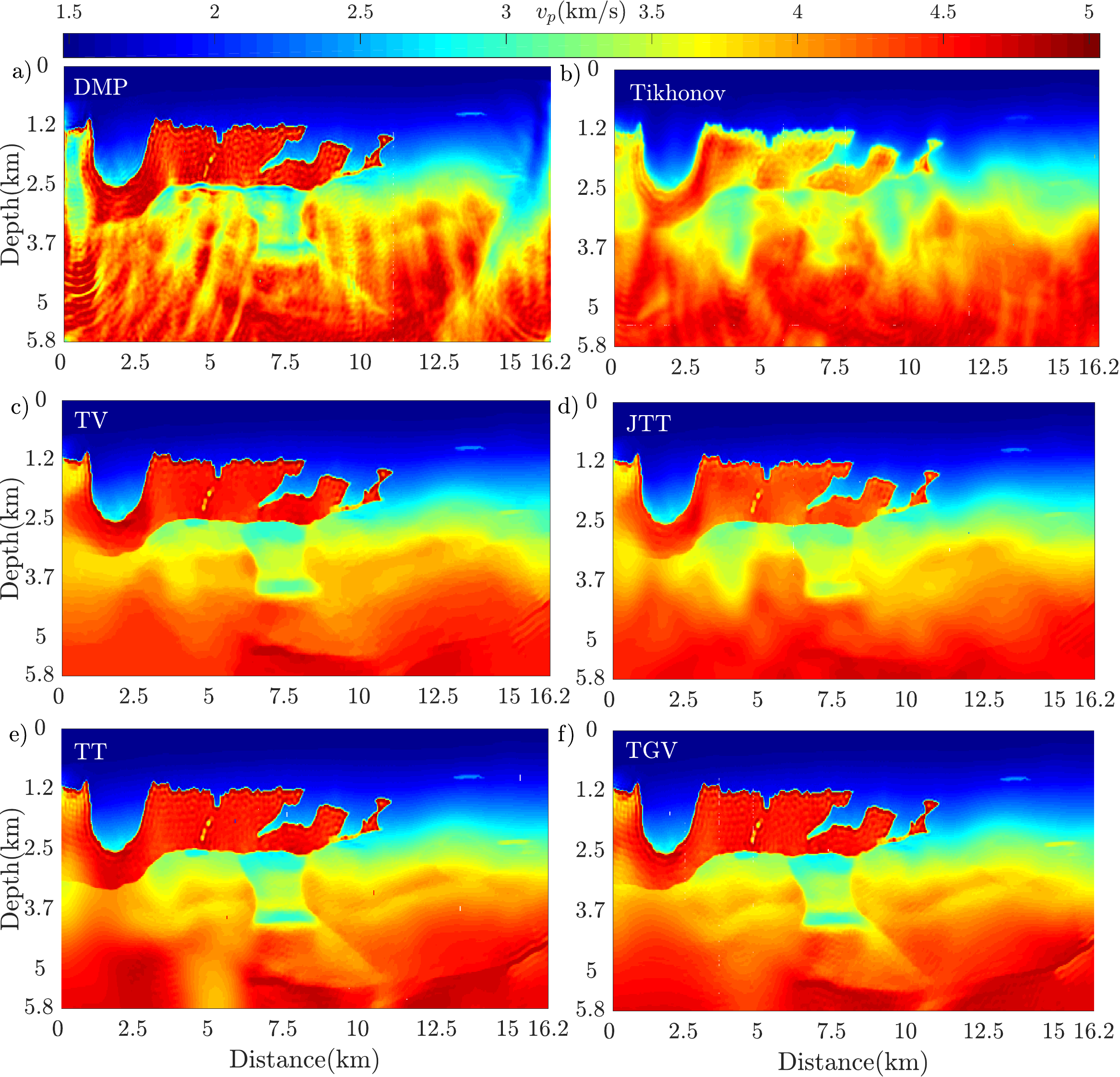}
\end{center}
\caption{2004 BP salt case study. Final inversion results with Fig. \ref{fig:bp_start}a-\ref{fig:bp_start}f as initial models. The panels are same as Fig. \ref{fig:bp_start}.}
\label{fig:bp}
\end{figure}

\begin{figure}[!h]
\begin{center}
\includegraphics[scale=0.45]{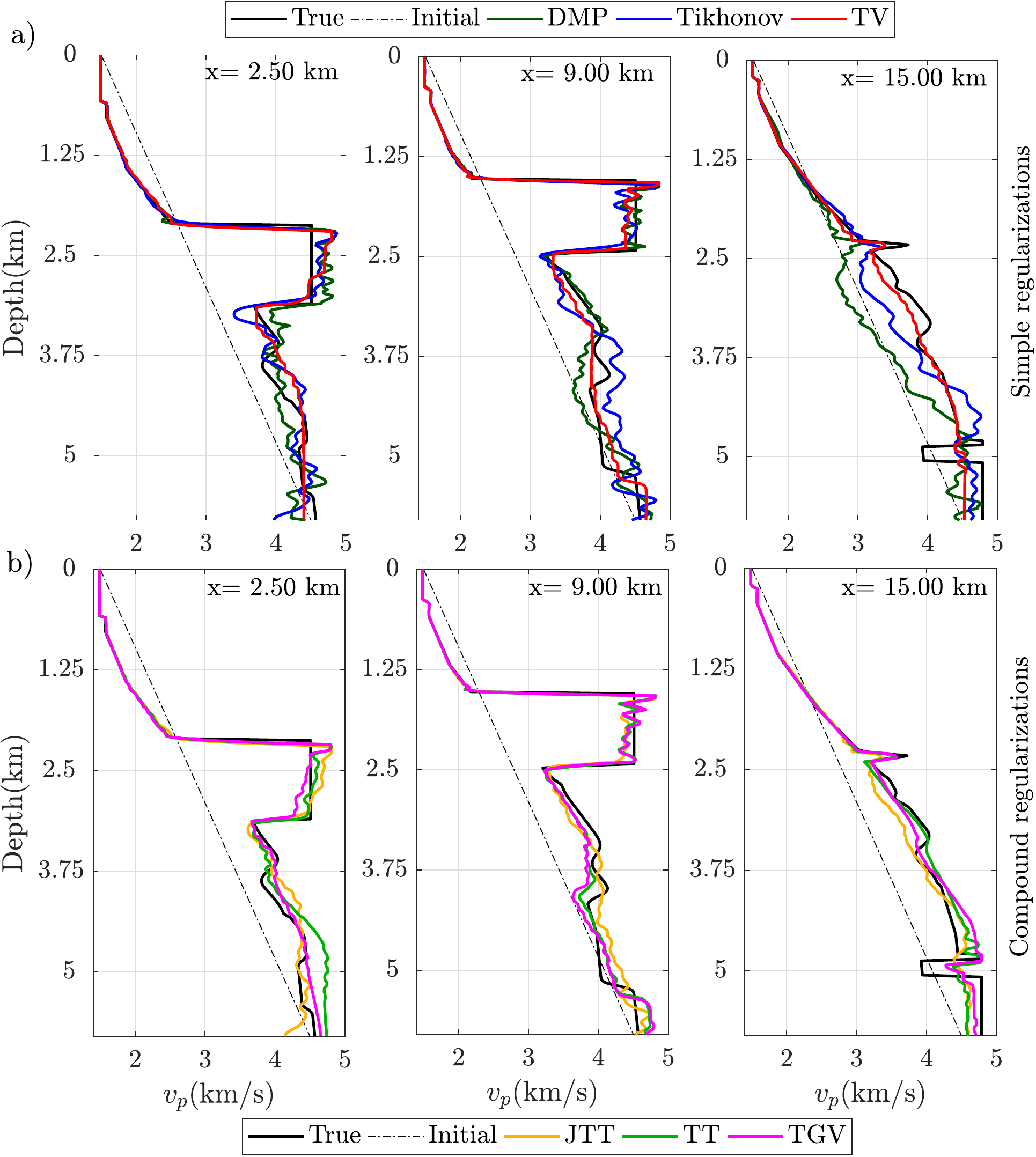}
\end{center}
\vspace{-2pt}
\caption{2004 BP salt case study. Direct comparison of final inversion results with Fig. \ref{fig:bp_start}a-\ref{fig:bp_start}f as initial models. The panels are same as Fig. \ref{fig:bp_start_log} for final results of Fig. \ref{fig:bp}.}
 \label{fig:bp_log}
 \vspace{-7pt}
\end{figure}

\begin{figure} 
\centering
   \includegraphics[width=0.48\textwidth]{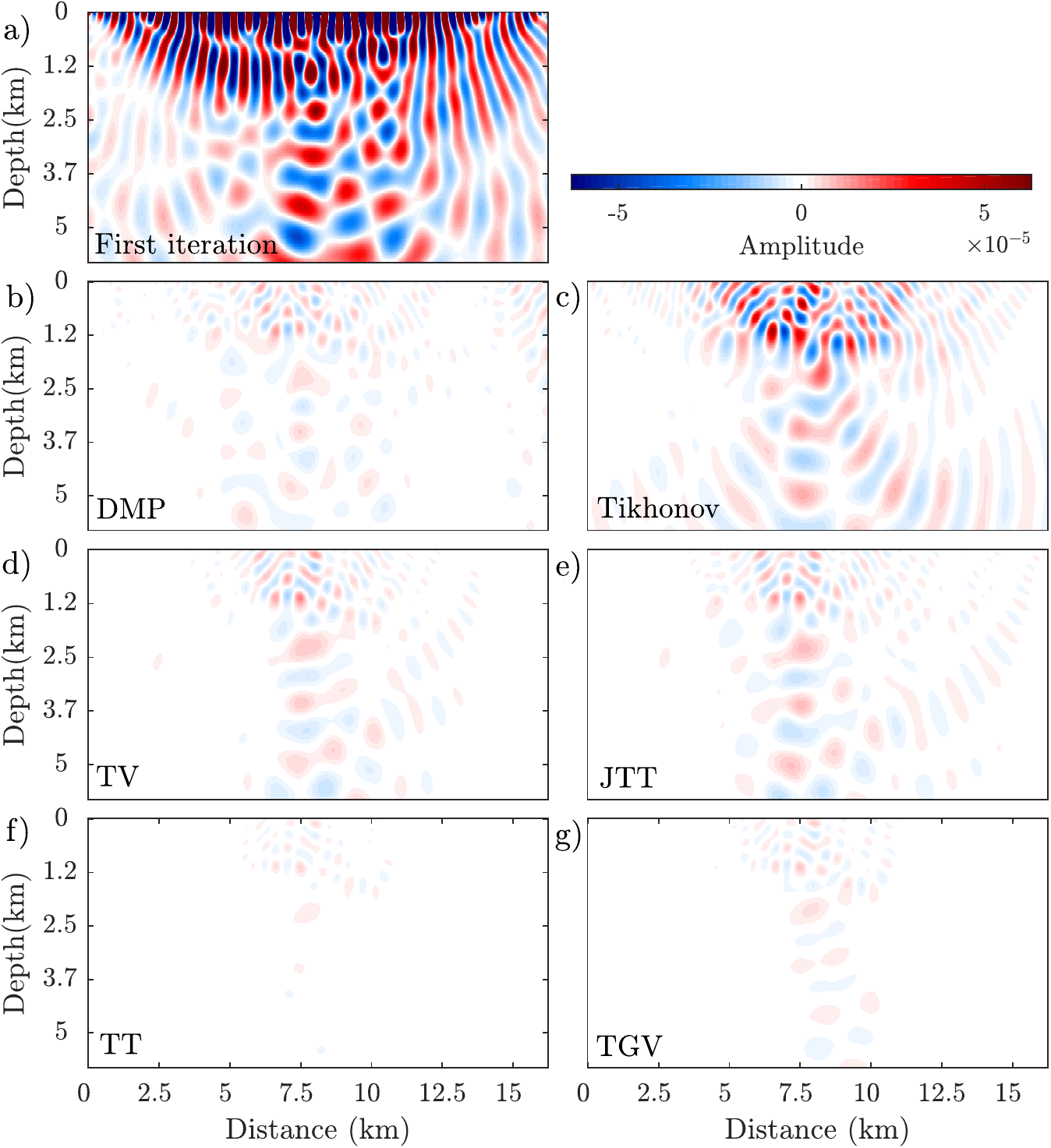} 
\caption{2004 BP salt case study. For a source at x=8.12~km and the 3-Hz frequency: (a) Real part of  wave-equation residual at first iteration of wavefield reconstruction, namely ($\bold{A}(\bold{m}^0)\bold{u}^1-\bold{b}$). (b-g) Real part of  wave-equation residual at the final iteration achieved respectively with DMP, Tikhonov, TV, JTT, TT, and TGV regularization.}
\label{fig:bp_res}
\end{figure}

\begin{figure} 
\includegraphics[width=0.48\textwidth]{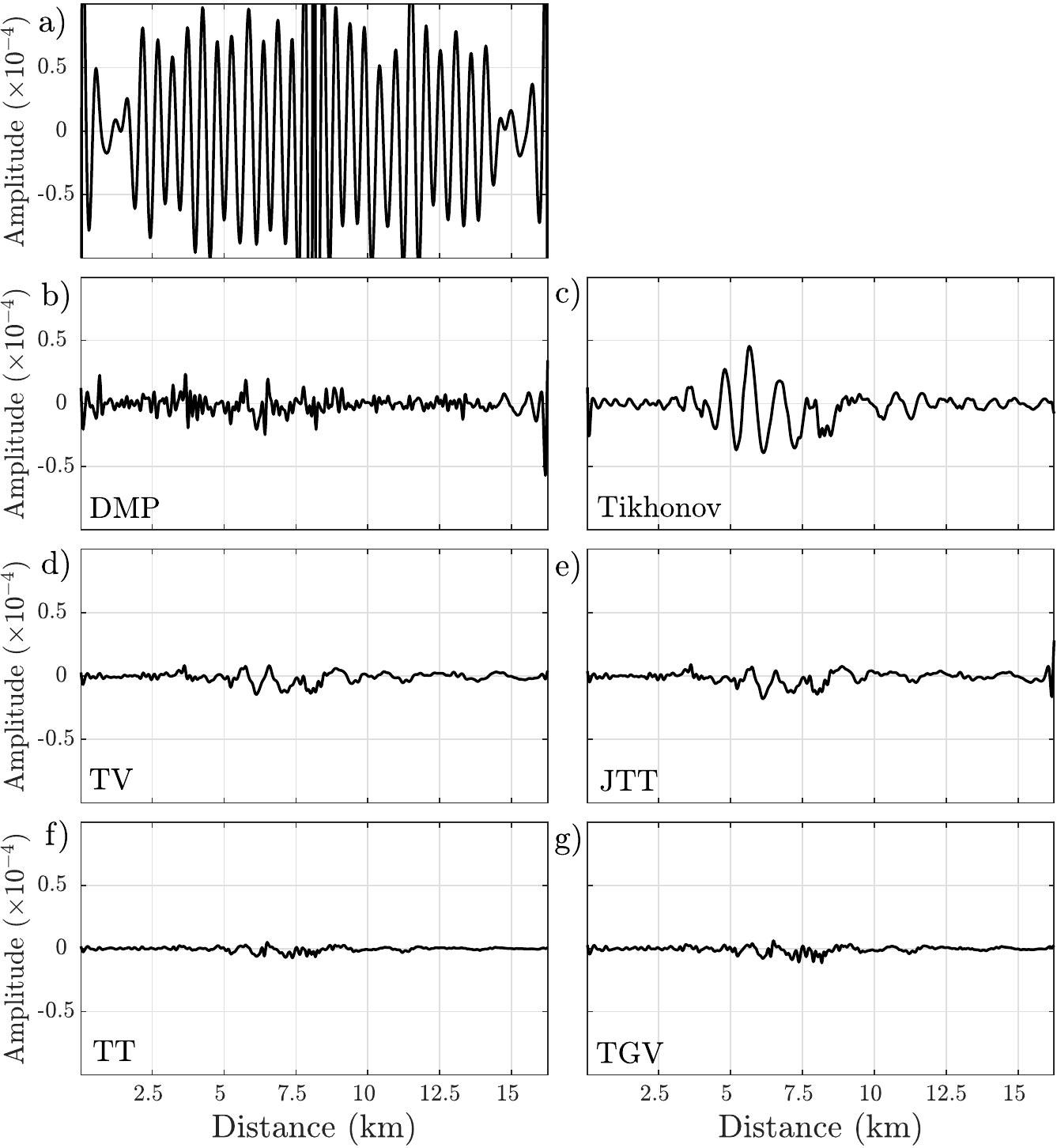}
\caption{2004 BP salt case study. For a source at x=8.12~km and the 3-Hz frequency: (a) Real part of data residual ($\bold{Pu}-\bold{d}$) at first  iteration. (b-g) Real part of data residual at the final iteration achieved respectively with DMP, Tikhonov, TV, JTT, TT, and TGV regularization.} 
\label{fig:bp_res1}
\end{figure}

We now consider a more realistic application with a target of the challenging 2004 BP salt model \cite{Billette_2004_BPB}. The 2004 BP salt model is representative of the geology of the deep offshore Gulf of Mexico and mainly consists of a simple background with a complex rugose multi-valued salt body, sub-salt slow velocity anomalies related to over-pressure zones and a fast velocity anomaly to the right of the salt body. 
The selected subsurface model is 16250~m wide and 5825~m deep, and is discretized with a 25~m grid interval (Fig. \ref{fig:bp_true}a). We used 108 sources spaced 150~m apart on the top side of the model. We perform forward modeling with a staggered-grid 9-point finite-difference method with PML boundary conditions using a 10~Hz Ricker wavelet as source signature. A line of receivers with a 25~m spacing are deployed at the surface leading to a stationary-receiver acquisition. 
We used small batches of two frequencies with one frequency overlap between two consecutive batches, moving from the low frequencies to the higher ones according to a classical frequency continuation strategy. The starting and final frequencies are 3~Hz and 13~Hz and the sampling interval in one batch is 0.5~Hz. The initial velocity model is a crude laterally-homogeneous velocity-gradient model with velocities ranging between 1.5 to 4.5~km/s, meaning that we initiate the inversion from scratch (Fig. \ref{fig:bp_true}b).
We start with inverting the first batch of frequencies (\{3, 3.5\}~Hz) with noiseless data using a maximum number of iteration equal to 45 as stopping criterion of iteration. To highlight the specific role of bound constraints, we activate them after 20 iterations.
To emphasize the effect of regularization, the result of bound constrained IR-WRI with a simple DMP regularization is shown in Fig. \ref{fig:bp_start}a, while the bound-constrained IR-WRI results with Tikhonov and TV regularizations are shown in Figs. \ref{fig:bp_start}b and \ref{fig:bp_start}c, respectively. Although the TV reconstruction is better than the Tikhonov one, it provides a velocity model which is far from the optimal one. A direct comparison between the true model, the starting model and the estimated models are shown in Fig.~\ref{fig:bp_start_log}a along three vertical logs at 2.5, 9.0 and 15.0~km distance (as depicted with dashed white lines in Fig. \ref{fig:bp_true}a).
We continue with compound regularization results which are shown in Fig. \ref{fig:bp_start}(d-f) and Fig.~\ref{fig:bp_start_log}b. Clearly, the TT regularizer better captures the long wavelengths of the salt body and the smooth subsalt background model.  The joint evolution in iterations of the observation-equation ($\|\bold{Pu-d}\|_2$) and wave-equation errors ($\|\bold{A(m)u-b}\|_2$) (Fig.~\ref{fig:bp_start_res}a) and the relative  model errors in iterations (Fig.~\ref{fig:bp_start_res}b) further confirm the relative performances of each regularizer during the inversion of the first frequency batch. Note the complex zigzag path followed by the inversion to jointly minimize the data residuals and the wave equation error in Fig.~\ref{fig:bp_start_res}(a-b). As already highlighted by \cite{Aghamiry_2019_IWR}, this results from the dynamic balancing in iterations of the observation-equation and wave-equation constraints performed by the dual updates with the data and source residuals.  \\
%
%
%
%
%
%
%
%
%
%
We continue the inversion at higher frequencies using the final models of the \{3, 3.5\}~Hz inversion, Fig.~\ref{fig:bp_start}a-f, as initial models when the stopping criteria is either $k_{max}=15$ or
\begin{eqnarray} \label{Stop}
\| \bold{A(m}^{k+1})\bold{u}^{k+1}-\bold{b}\|_2 \leq \varepsilon_b, ~~~~
 \|\bold{Pu}^{k+1}-\bold{d}\|_2 \leq \varepsilon_d,
\end{eqnarray} 
where $k_{max}$ denotes the maximum iteration count, $\varepsilon_b$=1e-3, and $\varepsilon_d$=1e-5. We perform three paths through the frequency batches to improve the IR-WRI results, using the final model of one path as the initial model of the next one (these cycles can be viewed as outer iterations of IR-WRI). The starting and finishing frequencies of the paths are [3.5, 6], [4, 8.5], [6, 13]~Hz respectively, where the first element of each pair shows the starting frequency and the second one is the finishing frequency.
The bound-constrained IR-WRI models obtained from noiseless data are shown in Fig.~\ref{fig:bp}. The number of iterations that have been performed with DMP regularization (Fig. \ref{fig:bp}a) is 426. Also, the number of iterations with Tikhonov (Fig. \ref{fig:bp}b) and TV regularizations (Fig. \ref{fig:bp}c) are 448 and 399, respectively, while it is 415 for the JTT regularization (Fig. \ref{fig:bp}d). Also, the number of iterations that have been performed with TT (Fig. \ref{fig:bp}e) and TGV (Fig. \ref{fig:bp}f) are 361 and 394, respectively. As for the inversion of the first batch, direct comparison between the true model, the starting model and the estimated models are shown in Fig.~\ref{fig:bp_log} along three vertical logs  at 2.5~km, 9.0~km and 15~km distance (vertical dashed lines in Fig.~\ref{fig:bp_true}a). The TT and TGV regularizers lead to high-quality velocity models, that capture both the fine-scale structure of the rugose large-contrast salt body and the high-velocity shallow anomaly on the right, as well as the smoother sub-salt background model including the low-velocity over-pressure structure. Moreover, the TT has a better convergence rate compare to the TGV. It is also worth noting the significant differences between the JTT and TT IR-WRI models in particular in the deep part of the model. \\
As a final quality control of the different IR-WRI models, it is instructive to check the wave-equation and data residuals left by the different regularization methods for the starting 3-Hz frequency (Fig. \ref{fig:bp_res} and \ref{fig:bp_res1}). The real part of wave-equation error (Fig \ref{fig:bp_res}) and data residuals (Fig \ref{fig:bp_res1}) are plotted at the first and final iterations of the inversion. Both the final data and source residuals suggest that the TT regularizer slightly outperforms the TGV counterpart at low frequencies.
 
%
%
%

\section{Conclusions} 
In this study, we first show how to implement efficiently different kinds of regularization and bound constraints in the wavefield reconstruction inversion method with the alternating direction method of multipliers (ADMM). Then, we show the capability of IR-WRI when equipped with compound Tikhonov and TV regularizations to reconstruct accurately large-contrast subsurface media when starting from very crude initial models. This compound regularization is suitable for seismic imaging of the subsurface as this later can be often represented by piecewise smooth media. We show that the infimal convolution (IC) of the Tikhonov and TV regularizers captures much more efficiently the blocky and smooth components of the subsurface than the convex combination of the two regularizers. It also outperforms the Tikhonov and TV regularizers when used alone. 
We also show how the infimal-convolution regularizer can be efficiently implemented by jointly updating the smooth and blocky subsurface components through variable projection. Alternatively, TGV regularized IR-WRI can be a suitable tool to reconstruct piecewise linear media and provides similar results than TT IR-WRI. We conclude that such hybrid regularizations in the extended search-space IR-WRI potentially provide a suitable framework to reconstruct without cycle skipping large-contrast subsurface media from ultra-long offset seismic data. It should also find applications in other fields of imaging sciences such as medical imaging.
\section{Acknowledgments}  
This study was partially funded by the SEISCOPE consortium (\textit{http://seiscope2.osug.fr}), sponsored by AKERBP, CGG, CHEVRON, EQUINOR, EXXON-MOBIL, JGI, PETROBRAS, SCHLUMBERGER, SHELL, SINOPEC and TOTAL. This study was granted access to the HPC resources of SIGAMM infrastructure (http://crimson.oca.eu), hosted by Observatoire de la C\^ote d'Azur and which is supported by the Provence-Alpes C\^ote d'Azur region, and the HPC resources of CINES/IDRIS/TGCC under the allocation A0050410596 made by GENCI."
%
%


%
\IEEEpeerreviewmaketitle

\bibliographystyle{IEEEtran}


\end{document}